\documentclass{article}
\usepackage{amsthm,amsmath,amssymb}
\usepackage{color}
\hoffset - 20 mm

\definecolor{dred}{rgb}{.8,0,0}
\definecolor{ddmagenta}{rgb}{0.7,0,0.9}
\definecolor{ddcyan}{rgb}{0,0.2,1.0}
\definecolor{dblue}{rgb}{0,0,0.7}

\newtheorem{teor}{Theorem}[section]
\newtheorem{defin}[teor]{Definition}
\newtheorem{lemm}[teor]{Lemma}
\newtheorem{osse}[teor]{Remark}
\newtheorem{prop}[teor]{Proposition}
\newtheorem{defi}[teor]{Definition}
\newtheorem{coro}[teor]{Corollary}
\newtheorem{prob}[teor]{Problem}

\newcommand{\bele}{\begin{lemm}\begin{sl}}
\newcommand{\enle}{\end{sl}\end{lemm}}
\newcommand{\bedef}{\begin{defi}\begin{sl}}
\newcommand{\eddef}{\end{sl}\end{defi}}
\newcommand{\bete}{\begin{teor}\begin{sl}}
\newcommand{\ente}{\end{sl}\end{teor}}
\newcommand{\beos}{\begin{osse}\begin{rm}}
\newcommand{\eddos}{\end{rm}\end{osse}}
\newcommand{\bepr}{\begin{prop}\begin{sl}}
\newcommand{\empr}{\end{sl}\end{prop}}
\newcommand{\bepro}{\begin{prob}\begin{rm}}
\newcommand{\empro}{\end{rm}\end{prob}}
\newcommand{\bede}{\begin{defin}\begin{sl}}
\newcommand{\edde}{\end{sl}\end{defin}}
\newcommand{\beco}{\begin{coro}\begin{sl}}
\newcommand{\enco}{\end{sl}\end{coro}}
\newcommand{\disp}{\displaystyle}

\newcommand{\quand}{\quad\text{and}\quad}

\newcommand{\quext}{\quad\text}

\newcommand{\de}{\partial}

\newcommand{\RR}{\mathbb{R}}
\newcommand{\EE}{\mathbb{E}}

\newcommand{\beeq}[1]{\begin{equation}\label{#1}}
\newcommand{\eddeq}{\end{equation}}

\newcommand{\beeqa}[1]{\begin{eqnarray}\label{#1}}
\newcommand{\eddeqa}{\end{eqnarray}}

\newcommand{\beal}[1]{\begin{align}\label{#1}}
\newcommand{\eddal}{\end{align}}

\newcommand{\bespl}[1]{\begin{split}\label{#1}}
\newcommand{\edspl}{\end{split}}

\newcommand{\bega}[1]{\begin{gather}\label{#1}}
\newcommand{\edga}{\end{gather}}

\newcommand{\beeqax}{\begin{eqnarray*}}
\newcommand{\eddeqax}{\end{eqnarray*}}

\def\qed{\ifmmode 
  \else \leavevmode\unskip\penalty9999 \hbox{}\nobreak\hfill
  \fi
  \quad\hbox{\hskip.5em\vrule width.4em height.6em depth.05em\hskip.1em}}
\def\endproofsym{\qed}
\newcommand{\dimbox}{\hbox{\hskip.5em\vrule width.4em height.6em depth.05em\hskip.1em}}
\renewenvironment{proof}[1][Proof]{\trivlist\item[\hskip\labelsep{\hskip0pt
    {\normalfont\scshape#1.}\hskip .321429\parindent}]\ignorespaces}
{\endproofsym\endtrivlist}
\def\endnobox{\def\endproofsym{}\end{proof}\def\endproofsym{\qed}}

\newcommand{\no}{\nonumber}

\newcommand{\beeqao}{\begin{eqnarray}\no}
\newcommand{\bealo}{\begin{align}\no}
\newcommand{\besplo}{\begin{split}\no}
\newcommand{\begao}{\begin{gather}\no}

\newcommand{\duav}[1]{\langle{#1}\rangle}

\def \no#1#2#3 {{\bf #1} (#3), #2.}
\def \eds#1#2#3 {#1, #2, #3.}

\newcommand{\dt}{\hbox{d}t}

\newcommand{\perogni}{\forall\,}
\newcommand{\esiste}{\exists\,}
\newcommand{\itt}{\int_0^t}

\newcommand{\io}{\int_\Omega}
\newcommand{\ibaro}{\int_{\overline\Omega}}

\newcommand{\iga}{\int_\Gamma}

\newcommand{\epsi}{\varepsilon}

\newcommand{\OO}{_{\Omega}}

\newcommand{\dn}{\partial_{\nu}}

\newcommand{\lhs}{left-hand side}
\newcommand{\rhs}{right-hand side}

\DeclareMathOperator{\deriv}{d}

\DeclareMathOperator{\dist}{dist}
\DeclareMathOperator{\dom}{dom}

\DeclareMathOperator{\sign}{sign}

\newcommand{\LDH}{L^2(0,T;H)}

\let\TeXchi\chi
\def\chi{{\setbox0 \hbox{\mathsurround0pt
$\TeXchi$}\hbox{\raise\dp0 \copy0 }}}

\newcommand{\ugat}{u_{\Gamma,t}}

\newcommand{\teta}{\vartheta}

\newcommand{\calH}{{\mathcal H}}

\newcommand{\calX}{{\mathcal X}}

\newcommand{\calA}{{\mathcal A}}

\newcommand{\calE}{{\mathcal E}}

\newcommand{\diX}{\dist_{\mathcal X}}
\newcommand{\AAA}{{\mathbf A}}

\newcommand{\calV}{{\mathcal V}}

\newcommand{\calW}{{\mathcal W}}

\newcommand{\chibar}{\overline{\chi}}
\newcommand{\ubar}{\overline{u}}

\newcommand{\Ga}{_\Gamma}
\newcommand{\barO}{\overline{\Omega}}

\newcommand{\ditau}{\deriv\!\tau}

\newcommand{\dit}{\deriv\!t}
\newcommand{\dis}{\deriv\!s}
\newcommand{\dix}{\deriv\!x}
\newcommand{\dx}{\deriv\!x}
\newcommand{\ds}{\deriv\!s}
\newcommand{\dm}{\deriv\!m}
\newcommand{\dS}{\deriv\!S}

\newcommand{\ddt}{\frac{\deriv\!{}}{\dit}}

\newcommand{\ii}{_\infty}
\newcommand{\iim}{_{\infty,\mu}}
\newcommand{\iig}{_{\infty,\Gamma}}

\newcommand{\tetaga}{{\teta}_{\Gamma}}
\newcommand{\chiga}{{\chi}_{\Gamma}}

\newenvironment{rcomm}{\color{ddmagenta} \textsf{R:}\,}{\color{ddmagenta}}
\newenvironment{secomm}{\color{dblue} \textsf{SE:}\,}{\color{dblue}}

\newenvironment{new}{\color{dred}}{\color{black}}
\newcommand{\bne}{\begin{new}}
\newcommand{\ene}{\end{new}}
\newcommand{\beroc}{\begin{rcomm}}
\newcommand{\eroc}{\end{rcomm}}

\newcommand{\bseg}{\begin{secomm}}
\newcommand{\eseg}{\end{secomm}}

\begin{document}



\title{The Penrose-Fife phase-field model\\
with coupled dynamic boundary conditions}

\author{Alain Miranville\thanks{The work of A.M. was partially supported 
  by the FP7-IDEAS-ERC-StG \#256872 (EntroPhase). Part of the work 
  has been developed while A.M.~was visiting the University of Milano
  from February 27 to March 4, 2012.}\\
Laboratoire de Math\'ematiques et Applications,\\
UMR CNRS 7348, Universit\'e de Poitiers - SP2MI,\\
Boulevard Marie et Pierre Curie,\\ 
F-86962 Chasseneuil Futuroscope Cedex, France\\
E-mail: {\tt Alain.Miranville@math.univ-poitiers.fr}\\
\and
Elisabetta Rocca\thanks{The work of E.R. was supported
   by the FP7-IDEAS-ERC-StG \#256872 (EntroPhase).}\\
Dipartimento di Matematica, Universit\`a di Milano,\\
Via Saldini 50, 20133 Milano, Italy\\
E-mail: {\tt elisabetta.rocca@unimi.it}\\
\and
Giulio Schimperna\thanks{The work of G.S.~was supported by the MIUR-PRIN Grant     
   2008ZKHAHN ``Phase transitions, hysteresis and multiscaling''
   and by the FP7-IDEAS-ERC-StG \#256872 (EntroPhase).}\\
Dipartimento di Matematica, Universit\`a di Pavia,\\
Via Ferrata~1, I-27100 Pavia, Italy\\
E-mail: {\tt giusch04@unipv.it}\\
\and
Antonio Segatti\thanks{The work of A.S.~was supported by the MIUR-PRIN Grant
   2008ZKHAHN ``Phase transitions, hysteresis and multiscaling''
   and by the FP7-IDEAS-ERC-StG \#256872 (EntroPhase).}\\
Dipartimento di Matematica, Universit\`a di Pavia,\\
Via Ferrata~1, I-27100 Pavia, Italy\\
E-mail: {\tt antonio.segatti@unipv.it}
}

\maketitle

\begin{abstract}
In this paper we derive, starting from the basic principles of
thermodynamics, an extended version of the nonconserved Penrose-Fife 
phase transition model, in which dynamic boundary conditions 
are considered in order to take into account interactions
with walls. Moreover, we study the well-posedness
and the asymptotic behavior
of the Cauchy problem for the PDE system associated to the model,
allowing the phase configuration of the material to be described
by a singular function.
\end{abstract}

\vspace{.4cm}


\vspace{4mm}

\noindent
{\bf AMS (MOS) subject clas\-si\-fi\-ca\-tion:} 
35K61, 35D30, 34B16, 74H40, 34K21, 80A22.
%
%


\section{Introduction}
\label{sec:intro}

In this paper we derive a model for phase transitions of Penrose-Fife type 
settled in a bounded domain $\Omega\subset \RR^3$.  
The peculiarity of our approach consists in the fact that we take into account the
relations between $\Omega$ and its exterior, including the effects of the 
interactions with the boundary into the free energy and
entropy functionals $\Psi$ and $S$
(cf.~\eqref{15}-\eqref{16} below) 
which drive the evolution of the system. A detailed derivation
of the model is carried out in Section \ref{sec:model}. The resulting PDEs 
system couples four nonlinear and singular evolution equations: two 
for the absolute temperature $\teta$ (one in the
bulk~$\Omega$ and the other on the boundary $\Gamma$) and two for the phase 
parameter $\chi$, which represents the local proportion of one
of the two phases:
\begin{equation}\label{eqtetaOmega}
\frac{\partial \teta}{\partial t}-\Delta {\left(-\frac1{\teta}\right) }
 +\lambda '_b(\chi ) \frac{\partial \chi }{\partial t}=h_b\ {\rm in}\ \Omega ,
\end{equation}
\begin{equation}\label{eqtetaGamma}
\frac{\partial \teta}{\partial t} -\Delta _\Gamma {\left(-\frac1{\teta}\right) }
 +\lambda '_\Gamma (\chi )\frac{\partial \chi }{\partial t}
  + \dn \left(-\frac1{\teta} \right)=h_\Gamma \ {\rm on}\ \Gamma ,
\end{equation}
\begin{equation}\label{eqchiOmega}
 \frac{\partial \chi}{\partial t}
  -\Delta \chi -s'_{0,b}(\chi )
   = - \frac{\lambda '_b(\chi )}{\teta }\ {\rm in}\ \Omega ,
\end{equation}
\begin{equation}\label{eqchiGamma}
\frac{\partial \chi}{\partial t}
-\Delta _\Gamma \chi -s'_{0,\Gamma }(\chi )
  + \dn\chi = - \frac{\lambda '_\Gamma (\chi )}{\teta }\ 
  {\rm on}\ \Gamma.
\end{equation}
Here, $\Delta$ stands for the Laplacian with respect to the space variables in $\Omega$, 
$\Delta_\Gamma$ denotes the Laplace-Beltrami operator on $\Gamma$,
and $\nu$ is the unit outer normal vector to $\Gamma$.
Moreover, $\lambda_b$ and $\lambda_\Gamma$ are two quadratic functions of 
$\chi$ related to the latent heat of the process,
$h_b$ and $h_\Gamma$ are two heat sources, respectively in the bulk and on the boundary, 
and $-s_{0,b}'$ and $-s_{0,\Gamma}'$ are two
nonlinear functions whose antiderivatives $-s_{0,b}$ and $-s_{0,\Gamma}$ correspond
to the {\sl configuration potentials} of the phase variable.
We admit the case when the domains of $-s_{0,b}$ and $-s_{0,\Gamma}$ 
are bounded, with the purpose of excluding the unphysical values
of the variable $\chi$. In this case, we shall speak of {\sl singular potentials}.
A relevant example is the so-called logarithmic potential,
given by 
\begin{equation}\label{logpot}
  -s(r)=(1+r)\log(1+r)+(1-r)\log(1-r) - \frac\delta2 r^2,
   \ \delta\ge 0.
\end{equation}
With the choice \eqref{logpot}, $\pm1$ denote the pure states and the values 
$\chi\not\in[-1,1]$ are penalized by identically assigning the value $-\infty$
to $s$ outside $[-1,1]$. In what follows, 
the potentials $-s_{0,b}$ and $-s_{0,\Gamma}$ 
will be split into the sum of (dominating) monotone parts 
$f$, and $f_\Gamma$ respectively, 
and of quadratic perturbations (cf.~the last term in \eqref{logpot}).
The literature devoted to the mathematical features of
phase transition models endowed with singular potentials 
is rather wide, also in specific relation with the 
Penrose-Fife model. Among the various contributions we quote
\cite{MZ} (for the Cahn-Hilliard equation), 
\cite{GPS1} (for the Caginalp phase-field model), 
\cite{ssz12} (for the Penrose-Fife system), and the 
references in these papers.

A notable feature of system \eqref{eqtetaOmega}-\eqref{eqchiGamma} 
is the occurrence of {\sl dynamic boundary conditions}. This type of conditions
has been proposed in the literature in different contexts (for instance, in
the framework of the Allen-Cahn and Cahn-Hilliard models) with the aim of describing
the interactions between the interior of a domain and the walls
(cf., e.g.~\cite{FiMD1}, \cite{FiMD2}, \cite{FRDGMMR}, \cite{MZCH}, 
and references therein). In particular, the case of singular potentials in the
context of the Cahn-Hilliard evolution has been recently analyzed in 
\cite{GMS1}, \cite{GMS2} and \cite{MZCH},
whereas the Caginalp phase-field 
system (cf.~\cite{Ca}) with dynamic boundary conditions 
on the phase parameter has been considered in several 
papers (cf., e.g.~\cite{CGGM}, \cite{CGM}, 
\cite{CM}, \cite{GGbis}), while only recently it has been 
coupled with dynamic boundary conditions for 
both the phase parameter and the temperature. 
Such a problem has been considered in \cite{gms} and 
\cite{CoGM}, where the well-posedness and the asymptotic
behavior (in terms of attractors) of solutions have been 
studied, also in the case of singular potentials.

As far as the Penrose-Fife model (cf.~\cite{PenroseFife}) is concerned, 
a vast literature is devoted to the well-posedness 
(cf., e.g.~\cite{CLau}, \cite{HSZ}, \cite{Lau}) 
and to the long-time behavior of solutions
both in term of attractors (cf., e.g.~\cite{IK}, \cite{RS}, \cite{Sc09})
and of convergence of single trajectories to stationary states
(cf.~\cite{FS}). Most of these contributions deal with 
Robin boundary conditions for the temperature, i.e.,
\begin{equation}\label{boutetastandard}
 \frac{\partial (1/\teta) }{\partial \nu }
 = \xi (\teta -\teta _s)\ {\rm on}\ \Gamma ,
\end{equation}
and no-flux conditions for the phase parameter:
\begin{equation}\label{chinoflux}
   \frac{\partial \chi }{\partial \nu }=0\ {\rm on}\ \Gamma.
\end{equation}
Relation \eqref{boutetastandard} establishes that the heat flux through 
the boundary $\Gamma$ is proportional to the difference between the internal 
and external temperature $\teta_s>0$ by a positive coefficient $\xi$,
while \eqref{chinoflux} prescribes that the boundary has no influence 
on the phase change process. 
Finally, we would like to quote the recent paper 
\cite{IKN} where a phase separation model of Penrose-Fife type with Signorini 
boundary condition has been studied.

The main novelties of our contribution stand in the fact that we 
can consider dynamic boundary conditions for a coupled phase-field 
system which may display a singular character both in $\teta$ and 
in $\chi$ in view of the terms $1/\teta$ and $s_{0,\Gamma}'(\chi)$, 
$s_{0,b}'(\chi)$ in \eqref{eqtetaOmega}-\eqref{eqchiGamma}. 
Physically speaking, the main difference between the dynamic boundary 
conditions \eqref{eqtetaGamma}, \eqref{eqchiGamma}
and the ``standard'' ones \eqref{boutetastandard}-\eqref{chinoflux}
stands in the fact that, in the case of dynamic b.c.'s, the 
walls of the the container have a significant effect on 
the phase transition process. This may include the case of two
phase changing substances in contact, where one occupies a small layer
surrounding the other one and this layer is approximated by a surface
in the mathematical model (concentrated capacity phenomenon, 
see, e.g., \cite{SaV}, \cite{Sc99} for more details).

A rigorous derivation of system \eqref{eqtetaOmega}-\eqref{eqchiGamma} 
starting from the basic laws of Thermodynamics is carried out in
the first part of the paper. Namely, the equations are obtained
by combining the free energy balance with the entropy dissipation
inequality and imposing physically realistic constitutive expressions
for the energy and entropy functionals.

The second part of the paper is devoted to the mathematical analysis
of the system in the framework of weak solutions. 
In comparison with similar models, the main mathematical difficulty 
consists here in the coupling between
\begin{enumerate}
\item the singularity of the model, in particular due to the presence of
 \begin{enumerate}
 \item the term $1/\teta$ both in the bulk equations \eqref{eqtetaOmega} and \eqref{eqchiOmega}
  and in the boundary equations \eqref{eqtetaGamma} and~\eqref{eqchiGamma},
 \item the (possibly) singular functions $s_{0,b}'$ in the phase equation \eqref{eqchiOmega} in the bulk,
  and $s_{0,\Gamma}'$ in the phase equation \eqref{eqchiGamma} on the boundary,
 \end{enumerate}
\item and the occurrence of dynamic boundary conditions.
\end{enumerate}
Actually, as already noted in the case of the Caginalp model 
with dynamic boundary conditions (see \cite{gms}), if we are 
in presence of singular potentials, a particular care is 
needed. For instance, in order for equation \eqref{eqchiOmega}
to make sense a.e.~in the space-time domain, 
the monotone part $f$ of $-s_{0,b}'$ needs to be controlled
from above by the monotone part $f\Ga$ of $-s_{0,\Gamma}'$
(cf.~\eqref{compatibility} below). In particular, this happens when 
the potentials $-s_{0,b}'$ and $-s_{0,\Gamma}'$ have the same 
effective domain (e.g., the set $[-1,1]$, as in 
the case \eqref{logpot}) and $|s_{0,b}'|$ explodes at most 
as fast as $|s_{0,\Gamma}'|$ as $\chi$ approaches the boundary 
of this domain (the values $\pm1$ in the specific case).

It is worth noticing that, in the case of 
completely general singular potentials 
(i.e., without any compatibility condition assumed), a weak solution
is still expected to exist; however, the equations 
\eqref{eqchiOmega}-\eqref{eqchiGamma} ruling the behavior of 
$\chi$ need to be interpreted in a weaker sense 
either by means of duality arguments or 
of variational inequality techniques (cf.~\cite{gms}, \cite{MZCH} 
for more details). This issue will be addressed in 
a forthcoming paper, where we also plan 
to weaken the regularity assumptions on the initial temperature
$\teta_0$. Actually, in addition to the natural conditions 
represented by the finiteness of the initial energy and entropy,
we will assume here that $\teta_0$ has the $L^2$-regularity.
This assumption, albeit meaningful, is not required
in the physical derivation of the model and can be then
considered to be somehow artificial. In other similar 
contexts (cf., e.g.~\cite{ssz12}), the $L^2$-regularity
of the initial temperature has been shown to be not 
necessary and, for this reason, we will try to remove
it also for this model.

After proving well-posedness of system \eqref{eqtetaOmega}-\eqref{eqchiGamma},
we shall analyze further properties of solutions. Proceeding along the lines
of \cite{ssz12} where the case of standard b.c.'s is treated,
we shall prove, by using a Moser iteration argument, that 
the temperature is uniformly bounded from below for
strictly positive times. Moreover, if $\teta_0$ is slightly more
summable (cf.~\eqref{hyptetaLp} below), we also have a uniform
upper bound. This kind of behavior occurs commonly in
the framework of parabolic equations with very-fast diffusion terms
(see \cite{vaz2006} and \cite{bonf_vaz2010} for the Cauchy 
problem in the whole space and the recent contributions \cite{ssz12}
and \cite{ssz13} for the bounded domain case with Neumann and dynamic boundary
conditions, respectively). 
In particular, in dimension three the exponent $p=3$ 
appears to be critical in the sense that solutions starting
with initial conditions in $L^p$ with $p>3$ become $L^\infty$
for strictly positive times. On the other hand, for $p<3$ the situation
is drastically different as the self-similar solution
in $\mathbb{R}^3$ ($(\cdot)_+$ denoting the positive part)
\begin{equation}
\label{self-similar}
\Theta(x,t):=\frac{2(T-t)_{+}^{1/2}}{\vert x\vert}
\end{equation}
shows. The smoothing effect for $p=3$ is currently an
open problem (see \cite{vaz2006} for further details).  
Moreover, it is worth noting that, at least when no 
external source is present, the regularization estimates
are also uniform with respect to time and give rise to
additional regularity properties for both components of
the solution.

Taking advantage of the regularization estimates, in the last
part of the paper we finally investigate the long-time behavior of 
trajectories. Namely, we are able to prove that,
at least in case \eqref{hyptetaLp} holds, any weak
solution admits a non-empty $\omega$-limit set which
only contains stationary solutions of the system.

A more precise characterization of the long-time behavior
relies on the structure of the set of steady state solutions,
which requires some further explanation.
Actually, integrating \eqref{eqtetaOmega} in space, using 
the boundary condition \eqref{eqtetaGamma}, 
and assuming zero external source, one readily sees that
the value
\begin{equation}\label{medie}
  \mu = \io\big( \teta + \lambda_b(\chi) \big)\,\dx
   + \iga\big( \teta + \lambda\Ga(\chi) \big)\,\dS,
 \end{equation}
representing the ``total mass'' of the internal energy,
is conserved in the evolution of the system.
Consequently, any limit point of a given weak solution
has also to respect the constraint \eqref{medie},
with $\mu$ depending only on the initial datum.
However, when $\mu$ is small, we are not able to
exclude that the set of stationary states satisfying
\eqref{medie} might contain temperatures $\teta$ being
arbitrarily close to $0$ (note that the stationary
formulation of \eqref{eqtetaOmega}-\eqref{eqtetaGamma}
simply prescribes that $\teta$ is a constant function).
Consequently, the \rhs s of the steady state
system associated to \eqref{eqchiOmega}-\eqref{eqchiGamma} 
might be arbitrarily large, which considerably weakens
the regularity properties of the set of its solutions.

On the contrary, it is easy to show that this situation
cannot occur when either $\mu$ is large enough or 
$\lambda'$ satisfies a suitable sign condition
(cf.~\eqref{segnola} below). If either property holds,
the set of stationary states is bounded in a very strong norm,
{\sl independently}\/ of the magnitude of the
initial data (but depending on the value of $\mu$). 
This fact, together with the precompactness of 
solution trajectories in the natural energy space
and with the existence of a coercive
Liapounov functional (namely, the energy~$\Psi$),
implies that the system~\eqref{eqtetaOmega}-\eqref{eqchiGamma}
admits a smooth global attractor, which is the last result 
we prove.


We conclude by giving the plan of the paper: in the next Section~\ref{sec:model} 
we detail the physical derivation of the model from the basic thermodynamical principles. 
In Section~\ref{sec:main} we introduce our precise concept of solutions and
state our main mathematical results related to well-posedness
and regularization properties of weak solutions.
The proofs are given in the subsequent Section~\ref{sec-proof-exiuni}. 
Finally, the long-time behavior of solutions is separately analyzed in 
the last Section~\ref{sec:long}.



\section{Derivation of the model}
\label{sec:model}

The Penrose-Fife model is derived by considering the free energy
density $w$ and the entropy density $s$, assuming that these
quantities depend both on the order parameter $\chi $ and the
absolute temperature $\teta $ (as in the original paper by
Penrose and Fife \cite{PenroseFife}), 
and imposing that the basic laws of Thermodynamics are satisfied.
Assuming that we always have sufficient regularity of the involved variables, 
we can then write
\begin{equation}\label{1}
 w=e-\teta s\ ({\rm Gibbs'\ relation}),
\end{equation}
where
\begin{equation}\label{2}
 s=- \frac{\partial w}{\partial \teta },
\end{equation}
and the internal energy density $e$ is defined by
\begin{equation}\label{3}
  e={{\partial {(\frac{w}{\teta })}}\over {\partial {({1\over \teta })}}}.
\end{equation}
Moreover, we assume that the total free energy
has the expression
\begin{equation}\label{4}
  \Psi (\chi ,\nabla \chi ,\teta )
   =\int _\Omega \left({{\kappa \teta}\over 2}\vert \nabla \chi \vert ^2+w(\teta ,\chi )\right)\,\dix,
\end{equation}
whereas the entropy functional is given by
\begin{equation}\label{5}
   S (\chi ,\nabla \chi ,e)
     =\int _\Omega \left(-{{\kappa }\over 2}\vert \nabla \chi \vert ^2+s(e,\chi )\right)\, \dix,
\end{equation}
where $\Omega $ is the domain occupied by the system and $\kappa>0$ denotes 
an interfacial energy coefficient.

The evolution equations for $\chi $ and $\teta $ are then obtained by stating
the relations
\begin{equation}\label{6}
  {{\partial \chi }\over {\partial t}}=K^\star {{\delta S}\over {\delta \chi }},\ K^\star >0,
\end{equation}
\begin{equation}\label{7}
{{\partial e}\over {\partial t}}=K\Delta {{\delta S}\over {\delta e}}+h,\ K>0,
\end{equation}
where $\delta $ denotes a variational derivative and $h$ is a source term.

Taking finally
\begin{equation}\label{8}
  e=c_0\teta +\lambda (\chi ),
\end{equation}
where $c_0>0$ stands for the specific heat of the system (assumed constant)
and $\lambda$ is the latent heat density, typically (cf.~\cite[p.~53]{PenroseFife})
given by
\begin{equation}\label{9}
  \lambda (r)=-ar^2+br+c,\ a>0,
\end{equation}
we find (see, e.g., \cite{MiranvilleRavello}; see also below)
\begin{equation}\label{10}
s(e,\chi )=c_0\ln \teta +s_0(\chi )+c_1.
\end{equation}
Here, $-s_0$ denotes a configuration potential, typically having
a double-well character (one can also consider a logarithmic 
double-well potential of the form \eqref{logpot}),
and such that $s''_0\le \delta$, $\delta\ge 0$, and
$c_1$ is a constant. These choices give rise to the Penrose-Fife
system:
\begin{equation}\label{11}
{{\partial \chi }\over {\partial t}}=K^\star \left(\kappa \Delta \chi +s'_0(\chi )
-{{c_0\lambda '(\chi )}\over \teta }\right),
\end{equation}
\begin{equation}\label{12}
c_0{{\partial \teta }\over {\partial t}}=-K\Delta {c_0\over \teta }-\lambda '(\chi ){{\partial \chi }\over {\partial t}}
+h.
\end{equation}
These equations are usually endowed with the boundary conditions
(cf.~the Introduction)
\begin{equation}\label{13}
{{\partial \chi }\over {\partial \nu }}=0\ {\rm on}\ \Gamma ,
\end{equation}
\begin{equation}\label{14}
{{\partial \teta }\over {\partial \nu }}=-\xi (\teta -\teta _s)\ {\rm on}\ \Gamma ,\ \xi ,\ \teta _s>0,
\end{equation}
where $\Gamma =\partial \Omega $ and $\nu $ is the unit outer normal vector to $\Gamma $.

Now, in order to take into account the interactions with the exterior of $\Omega$ 
(e.g., the walls), it is natural, following \cite{KEMRSBD} (see also \cite{CoGM} 
and \cite{GoMS}), to add a boundary contribution to the total free energy and
to take, in place of \eqref{4},
\begin{equation}\label{15}
  \Psi =\Psi (\chi ,\nabla \chi ,\nabla _\Gamma \chi , \teta )
   =\int _\Omega \left({{\kappa _b\teta}\over 2}\vert \nabla \chi \vert ^2+w_b(\teta ,\chi )\right)\, \dix
\end{equation}
\begin{equation*}
 + \int _\Gamma \left({{\kappa _\Gamma \teta}\over 2}\vert \nabla _\Gamma \chi \vert ^2
   + w_\Gamma (\teta ,\chi )\right)\, \dS ,\ \kappa _b,\ \kappa _\Gamma >0,
\end{equation*}
where $\nabla _\Gamma $ is the surface gradient and $w_b$ and $w_\Gamma $ 
are the bulk and surface free energy densities, respectively. Similarly, 
it is reasonable, in view of \eqref{15}, to introduce the generalized
entropy functional
\begin{equation}\label{16}
  S=S(\chi ,\nabla \chi ,\nabla _\Gamma \chi , e)
   =\int _\Omega \left(-{{\kappa _b}\over 2}\vert \nabla \chi \vert ^2+s_b(e,\chi )\right)\, \dix
   +\int _\Gamma \left(-{{\kappa _\Gamma }\over 2}\vert \nabla _\Gamma \chi \vert ^2+s_\Gamma (e,\chi )\right)\, \dS,
\end{equation}
%
%
where $s_b$ and $s_\Gamma $ are the bulk and surface entropy densities, respectively.

As before, we assume that Gibbs' relation holds, namely,
\begin{equation}\label{17}
  w_b=e-\teta s_b\ {\rm in}\ \Omega ,
\end{equation}
\begin{equation}\label{18}
  w_\Gamma =e-\teta s_\Gamma \ {\rm on}\ \Gamma ,
\end{equation}
and that
\begin{equation}\label{19}
  s_b=-{{\partial w_b}\over {\partial \teta }}\ {\rm in}\ \Omega ,
\end{equation}
\begin{equation}\label{20}
s_\Gamma =-{{\partial w_\Gamma }\over {\partial \teta }}\ {\rm on}\ \Gamma ,
\end{equation}
\begin{equation}\label{21}
e={{\partial {({w_b\over \teta })}}\over {\partial {({1\over \teta })}}}\ {\rm in}\ \Omega ,
\end{equation}
\begin{equation}\label{22}
e={{\partial {({w_\Gamma \over \teta })}}\over {\partial {({1\over \teta })}}}\ {\rm on}\ \Gamma.
\end{equation}
We now note that, in view of \eqref{16},
\begin{equation}\label{23}
  {{\delta S}\over {\delta \chi }}
   =\kappa _b\Delta \chi +{{\partial s_b}\over {\partial \chi }}\ {\rm in}\ \Omega ,
\end{equation}
\begin{equation}\label{24}
  {{\delta S}\over {\delta \chi }} 
  = \kappa _\Gamma \Delta _\Gamma \chi 
   +{{\partial s_\Gamma }\over {\partial \chi }}
    -\kappa _b{{\partial \chi }\over {\partial \nu }}\ {\rm on}\ \Gamma ,
\end{equation}
where $\Delta _\Gamma $ is the Laplace-Beltrami operator.

Then, assuming that, as in the classical model,
\begin{equation}\label{25}
{{\partial \chi }\over {\partial t}}=K^\star {{\delta S}\over {\delta \chi }},\
\ K^\star >0,
\end{equation}
we obtain the equations
\begin{equation}\label{26}
{{\partial \chi }\over {\partial t}}=K^\star \left(\kappa _b\Delta \chi +{{\partial s_b}\over {\partial \chi }}\right)\
{\rm in}\ \Omega ,
\end{equation}
\begin{equation}\label{27}
{{\partial \chi }\over {\partial t}}=K^\star \left(\kappa _\Gamma \Delta _\Gamma \chi +{{\partial s_\Gamma }\over {\partial \chi }}
-\kappa _b{{\partial \chi }\over {\partial \nu }}\right)\
{\rm on}\ \Gamma.
\end{equation}

Next, in order to describe the evolution of the
temperature, we generalize \eqref{7} as follows. We introduce, for
$U=\begin{pmatrix}u|_\Omega \\v|_\Gamma \end{pmatrix}$ 
(regular enough at this stage), the linear operator
$\mathbf A$ defined by
\begin{equation}\label{28}\displaystyle
\mathbf A U=\begin{pmatrix}
-\Delta u|_\Omega \\
-\Delta _\Gamma v|_\Gamma +{{\partial u|_\Omega }\over {\partial \nu }}|_\Gamma
\end{pmatrix}
\end{equation}
and write that
\begin{equation}\label{29}
{\partial \over {\partial t}}
\displaystyle
\begin{pmatrix}
e|_\Omega \\
e|_\Gamma
\end{pmatrix}
=-K\mathbf A
\begin{pmatrix}
{{\delta S}\over {\delta e}}|_\Omega \\
{{\delta S}\over {\delta e}}|_\Gamma
\end{pmatrix}
+\mathbf H,
\end{equation}
where $\mathbf H=\begin{pmatrix}h_b\\ h_\Gamma \end{pmatrix}$ is a forcing term,
with $h_b$ and $h_\Gamma$ standing for the bulk and boundary heat sources, 
respectively. Introducing the linear operator $\mathbf A$ is natural 
when considering dynamic boundary conditions (see \cite{CoGM} and \cite{GoMS}). Noting that
\begin{equation}\label{30}
 {{\delta S}\over {\delta e}}={{\partial s_b}\over {\partial e}}\ {\rm in}\ \Omega ,
\end{equation}
\begin{equation}\label{31}
{{\delta S}\over {\delta e}}={{\partial s_\Gamma }\over {\partial e}}\ {\rm on}\ \Gamma ,
\end{equation}
we deduce the equations
\begin{equation}\label{32}
{{\partial e}\over {\partial t}}=K\Delta {{\partial s_b}\over {\partial e}}+h_b\ {\rm in}\ \Omega ,
\end{equation}
\begin{equation}\label{33}
{{\partial e}\over {\partial t}}=K\Delta _\Gamma {{\partial s_\Gamma }\over {\partial e}}+h_\Gamma
-K{{\partial }\over {\partial \nu }}{{\partial s_b}\over {\partial e}}\ {\rm on}\ \Gamma.
\end{equation}

We then assume that (see \eqref{8})
\begin{equation}\label{34}
e=c_{0,b}\teta +\lambda _b(\chi )\ {\rm in}\ \Omega ,\ c_{0,b}>0,
\end{equation}

\begin{equation}\label{35}
e=c_{0,\Gamma }\teta +\lambda _\Gamma (\chi )\ {\rm on}\ \Gamma ,\ c_{0,\Gamma }>0.
\end{equation}

\noindent Noting that it follows from \eqref{21} and \eqref{22} that

\begin{equation}\label{36}
e=-\teta ^2{{\partial ({w_b\over \teta })}\over {\partial \teta }}\ {\rm in}\ \Omega ,
\end{equation}

\begin{equation}\label{37}
e=-\teta ^2{{\partial ({w_\Gamma \over \teta })}\over {\partial \teta }}\ {\rm on}\ \Gamma ,
\end{equation}

\noindent we have, freezing $\chi$ and integrating between $\teta _0$ and $\teta $, $\teta _0$, $\teta >0$,
\begin{equation}\label{38}
{{w_b(\teta, \chi)}\over \teta }-{{w_b(\teta _0, \chi)}\over \teta _0}=-\int _{\teta _0}^\teta
({c_{0,b}\over \tau }+{{\lambda _b(\chi )}\over \tau ^2})\, \ditau \ {\rm in}\ \Omega ,
\end{equation}
\begin{equation}\label{39}
{{w_\Gamma (\teta,\chi )}\over \teta }-{{w_\Gamma (\teta _0,\chi)}\over \teta _0}=-\int _{\teta _0}^\teta
({c_{0,\Gamma }\over \tau }+{{\lambda _\Gamma (\chi )}\over \tau ^2})\, \ditau \ {\rm on}\ \Gamma ,
\end{equation}
which we rewrite in the form
\begin{equation}\label{40}
w_b(\teta ,\chi )=-c_{0,b}\teta \ln \teta +c_{0,b}\teta \ln \teta _0+\lambda _b(\chi )
-\teta s_{0,b}(\chi )\ {\rm in}\ \Omega ,
\end{equation}
\begin{equation}\label{41}
w_\Gamma (\teta ,\chi )=-c_{0,\Gamma }\teta \ln \teta +c_{0,\Gamma }\teta \ln \teta _0+\lambda _\Gamma (\chi )
-\teta s_{0,\Gamma }(\chi )\ {\rm on}\ \Gamma ,
\end{equation}
where we have defined $s_{0,b}:=(w_b(\teta_0, \chi)-\lambda_b(\chi))/\teta_0$ and $s_{0,\Gamma}:=(w_\Gamma(\teta_0, \chi)-\lambda_\Gamma(\chi))/\teta_0$.

We finally deduce from \eqref{19}-\eqref{20} and \eqref{40}-\eqref{41} that

\begin{equation}\label{42}
s_b(e,\chi )=c_{0,b}\ln \teta +s_{0,b}(\chi )+c_{1,b},\ c_{1,b}=-c_{0,b}\ln \teta _0+c_{0,b}\ {\rm in}\ \Omega ,
\end{equation}

\begin{equation}\label{43}
s_\Gamma (e,\chi )=c_{0,\Gamma }\ln \teta +s_{0,\Gamma }(\chi )+c_{1,\Gamma },\ c_{1,\Gamma }
=-c_{0,\Gamma }\ln \teta _0+c_{0,\Gamma }\ {\rm on}\ \Gamma.
\end{equation}

We now note that it follows from \eqref{34}-\eqref{35} and \eqref{42}-\eqref{43} that

\begin{equation}\label{44}
s_b(e,\chi )=c_{0,b}\ln \Big({{e-\lambda _b(\chi )}\over {c_{0,b}}}\Big)+s_{0,b}(\chi )+c_{1,b}
\ {\rm in}\ \Omega ,
\end{equation}

\begin{equation}\label{45}
s_\Gamma (e,\chi )=c_{0,\Gamma }\ln \Big({{e-\lambda _\Gamma (\chi )}\over {c_{0,\Gamma }}}\Big)
+s_{0,\Gamma }(\chi )+c_{1,\Gamma }
\ {\rm on}\ \Gamma ,
\end{equation}

\noindent which yields

\begin{equation}\label{46}
{{\partial s_b}\over {\partial e}}={1\over \teta },\ {{\partial s_b}\over {\partial \chi }}=-
{{\lambda '_b(\chi )}\over \teta }+s'_{0,b}(\chi )\ {\rm in}\ \Omega ,
\end{equation}

\begin{equation}\label{47}
{{\partial s_\Gamma }\over {\partial e}}={1\over \teta },\ {{\partial s_\Gamma }\over {\partial \chi }}=-
{{\lambda '_\Gamma (\chi )}\over \teta }+s'_{0,\Gamma }(\chi )\ {\rm on}\ \Gamma.
\end{equation}

We finally deduce from \eqref{30}-\eqref{33} and \eqref{46}-\eqref{47} the following Penrose-Fife system
with dynamic boundary conditions:
\begin{equation}\label{48}
{{\partial \chi }\over {\partial t}}=K^\star \left(\kappa _b\Delta \chi +s'_{0,b}(\chi )-{{\lambda '_b(\chi )}
\over \teta }\right)\ {\rm in}\ \Omega ,
\end{equation}
\begin{equation}\label{49}
{{\partial \chi }\over {\partial t}}=K^\star \left(\kappa _\Gamma \Delta _\Gamma \chi +s'_{0,\Gamma }(\chi )-{{\lambda '_\Gamma (\chi )}
\over \teta }-\kappa _b{{\partial \chi }\over {\partial \nu }}\right)\ {\rm on}\ \Gamma ,
\end{equation}
\begin{equation}\label{50}
c_{0,b}{{\partial \teta }\over {\partial t}}=-K\Delta {1\over \teta }-\lambda '_b(\chi ){{\partial \chi }
\over {\partial t}}+h_b\ {\rm in}\ \Omega ,
\end{equation}
\begin{equation}\label{51}
c_{0,\Gamma }{{\partial \teta }\over {\partial t}}=-K\Delta _\Gamma {1\over \teta }-\lambda '_\Gamma (\chi ){{\partial \chi }
\over {\partial t}}+h_\Gamma +K{\partial \over {\partial \nu }}{1\over \teta }\ {\rm on}\ \Gamma.
\end{equation}

Taking, for simplicity, all constants equal to one, 
\eqref{48}-\eqref{51} reduces to
\begin{equation}\label{52}
  {{\partial \chi }\over {\partial t}}=\Delta \chi +s'_{0,b}(\chi )-{{\lambda '_b(\chi )}
   \over \teta }\ {\rm in}\ \Omega ,
\end{equation}
\begin{equation}\label{53}
{{\partial \chi }\over {\partial t}}=\Delta _\Gamma \chi +s'_{0,\Gamma }(\chi )-{{\lambda '_\Gamma (\chi )}
\over \teta }-{{\partial \chi }\over {\partial \nu }}\ {\rm on}\ \Gamma ,
\end{equation}
\begin{equation}\label{54}
{{\partial \teta }\over {\partial t}}=-\Delta {1\over \teta }-\lambda '_b(\chi ){{\partial \chi }
\over {\partial t}}+h_b\ {\rm in}\ \Omega ,
\end{equation}
\begin{equation}\label{55}
{{\partial \teta }\over {\partial t}}=-\Delta _\Gamma {1\over \teta }-\lambda '_\Gamma (\chi ){{\partial \chi }
\over {\partial t}}+h_\Gamma +{\partial \over {\partial \nu }}{1\over \teta }\ {\rm on}\ \Gamma ,
\end{equation}
which can also be rewritten in the following compact form:
\begin{equation}\label{56}
\displaystyle
{\partial \over {\partial t}}\begin{pmatrix}
\chi |_\Omega \\ \chi |_\Gamma \end{pmatrix}
=-\mathbf A \begin{pmatrix}
\chi |_\Omega \\ \chi |_\Gamma \end{pmatrix}
+\begin{pmatrix} s'_{0,b}(\chi )\\ s'_{0,\Gamma }(\chi ) \end{pmatrix}
-\begin{pmatrix} {{\lambda '_b(\chi )}
\over \teta }\\ {{\lambda '_\Gamma (\chi )}
\over \teta }\end{pmatrix},
\end{equation}
\begin{equation}\label{57}
\displaystyle
{\partial \over {\partial t}}\begin{pmatrix}
\teta |_\Omega \\ \teta |_\Gamma \end{pmatrix}
=\mathbf A \begin{pmatrix}
{1\over \teta }|_\Omega \\ {1\over \teta }|_\Gamma \end{pmatrix}
-\begin{pmatrix}\lambda '_b(\chi ){{\partial \chi }
\over {\partial t}}\\ \lambda '_\Gamma (\chi ){{\partial \chi }
\over {\partial t}}\end{pmatrix}+\mathbf H.
\end{equation}
The remainder of the paper is devoted to the mathematical analysis of system
\eqref{56}-\eqref{57} in the framework of weak solutions.



\section{Main assumptions and preliminary results}
\label{sec:main}

We introduce here our main assumptions, together
with several mathematical tools which are needed
in order to provide a precise analytical statement
of our results.

We let $\Omega$ be a sufficiently smooth, bounded, and connected
domain in $\RR^3$ with boundary $\Gamma$. 
We set $\barO:=\Omega \cup \Gamma$, 
set $H:=L^2(\Omega)$, and denote by
$(\cdot,\cdot)$ the scalar product both in $H$ and
in $H^3$ and by $\|\cdot\|$ the related norm.
Next, we set $V:=H^1(\Omega)$ and denote by
$V'$ the (topological) dual of $V$.
The duality between $V'$ and $V$
will be indicated by $\langle\cdot,\cdot\rangle$.
Identifying $H$ with $H'$ through the scalar product
of $H$, it is then well known that
$V\subset H\subset V'$ with continuous
and dense inclusions. In other words,
$(V,H,V')$ constitutes a Hilbert triplet
(see, e.g., \cite{Li}). Such a triplet
is usually used for stating
{\sl weak formulations}\/ of elliptic
or parabolic problems defined on $\Omega$.

However, since system \eqref{56}-\eqref{57}
also includes equations defined on $\Gamma$, we need
to introduce some further spaces taking also boundary
contributions into account. Thus, we set $H\Ga:=L^2(\Gamma)$,
$V\Ga:=H^1(\Gamma)$, and denote by
$(\cdot,\cdot)\Ga$ the scalar product in $H\Ga$,
by $\|\cdot\|\Ga$ the corresponding norm,
and by $\langle\cdot,\cdot\rangle\Ga$
the duality between $V'\Ga$ and $V\Ga$.
In general, the symbol $\|\cdot\|_{X}$
indicates the norm
in the generic (real) Banach space $X$ and
$\langle\cdot,\cdot\rangle_X$
stands for the duality between $X'$ and $X$.
We also denote by $\nabla\Ga$ the tangential gradient
on $\Gamma$ and by $\Delta\Ga$ the Laplace-Beltrami
operator. Thus, we can define the spaces
\begin{equation}\label{defiHVGa}
 \calH:=H\times H_\Gamma \quand
   \calV:=\big\{z\in V:~z|\Ga \in V\Ga\big\}.
\end{equation}
Here and in the following, $z|\Ga$, or also $z\Ga$, will denote
the trace of $z$ in the sense of a suitable trace theorem.
Next, we introduce the $\calH$-scalar product in the
following natural way:
\begin{equation}\label{defiproscaH}
  \big((k,\kappa),(s,\sigma)\big)_{\calH}:=
    (k,s)+(\kappa,\sigma)\Ga.
\end{equation}
Then, we set, respectively on $\calV$ and on $V$,
\begin{align}\label{defiproscaV}
  (\!(z,w)\!)_{\calV} & := \io (\nabla z\cdot\nabla w) \, \dx
   + \iga \big( z|\Ga w|\Ga
    + \nabla\Ga z|\Ga \cdot \nabla\Ga w|\Ga\big) \, \dS,\\
 \label{defiproscaV0}
  (\!(z,w)\!)_{V} & :=\io (\nabla z\cdot\nabla w) \, \dx
   + \iga z|\Ga w|\Ga \,\dS,
\end{align}
together with the related norms $\| \cdot \|_{\calV}$, $\| \cdot \|_{V}$.
It is not difficult to prove (see, e.g., \cite[Lemma~2.1]{MuRo}) 
that the space $\calV$ is dense in $\calH$.
Concerning the scalar product in \eqref{defiproscaV0},
we notice that the related norm $\|\cdot\|_V$ is
obviously equivalent to the usual one. 
We also set
\begin{equation}\label{defiWGa}
  \calW:=\big\{z\in \calV:~z\in H^2(\Omega),~
   z|\Ga\in H^2(\Gamma)\big\}
\end{equation}
and endow this space with the graph norm,
so that $\calW\subset \calV$, continuously and
compactly.

The above defined functional spaces allow
to introduce some elliptic operators. We set:
\begin{align}\label{defiA}
  A:V\to V', \qquad
   & \langle A z_1, z_2 \rangle
    :=\io \nabla z_1 \cdot \nabla z_2  \, \dx,\\
 \label{defiAGa}
  A\Ga:V\Ga\to V\Ga', \qquad
   & \langle A\Ga \zeta_1, \zeta_2 \rangle\Ga
    :=\iga \nabla\Ga \zeta_1 \cdot \nabla\Ga \zeta_2  \, \dS,\\
 \label{defiAcal}
  \calA:\calV \to \calV', \qquad
   & \langle \calA z_1, z_2 \rangle_\calV
    := \langle A z_1, z_2 \rangle
    + \langle A\Ga z_{1,\Gamma}, z_{2,\Gamma} \rangle\Ga.
\end{align}
We shall also use the
operator $\AAA$ defined in \eqref{28}. Note that $\AAA$ 
can be interpreted as an operator defined on $\calW$ and taking
values in $\calH$, thanks to the trace theorem
for normal derivatives.

In what follows, we will use the following convention:
as far as equations on $\Omega$ are concerned,
the elements of $\calV$ will be
interpreted as functions defined on $\Omega$ with
the proper regularity.
When, instead, as in most cases in the paper,
a system defined on $\Omega\times \Gamma$ is
considered, then the elements of $\calV$ will be
considered as {\sl pairs}\/ of functions $(z,z|\Ga)$.
In other words, $\calV$ will be identified
with a (closed) subspace of the product space
$H^1(\Omega)\times H^1(\Gamma)$. Analogously,
$V$ will be identified with a subspace of 
$H^1(\Omega)\times H^{1/2}(\Gamma)$,
in view of the trace theorem.
If we have, instead, $h\in\calH$, $h$ will be often thought
as a pair of functions belonging, respectively,
to $H$ and to $H\Ga$, and both denoted by the 
same letter $h$. Of course, if we do not have
additional regularity, the second component of $h$ needs
not be the trace of the first one.
Identifying $\calH$ with $\calH'$ 
{\sl through the scalar product \eqref{defiproscaH}}, 
we obtain the chain of continuous
and dense (thanks to the density of $\calV$ into $\calH$,
to the density of $H^2(\Omega)$ into $V$, and to
the continuous inclusion $H^2(\Omega)\subset \calV$)
embeddings
\begin{equation}\label{VVHVV}
  \calV\subset V\subset \calH\subset V'\subset \calV'.
\end{equation}
In particular, the relation
\begin{equation}\label{HVnorms}
  \big((k,\kappa),(z,z|\Ga)\big)_{\calH}
   = \io kz \,\dx + \iga \kappa z|\Ga \,\dS
   = \langle (k,\kappa),z\rangle_{\calV}
\end{equation}
holds for any $z\in \calV$ and $(k,\kappa)\in\calH$.
Of course, an analogous relation could be stated
for $z\in V$. 

%
Next, for any function, or functional $z$, defined on
$\Omega$, we set
\begin{equation}\label{defim}
  m\OO(z):=\frac1{|\Omega|}\io z \,\dx,
\end{equation}
where the integral is substituted with
the duality $\langle z,1 \rangle$ in case, e.g.,
$z\in V'$. We also define
the measure $\dm$ given by
\begin{equation}\label{defidm}
  \ibaro v \dm:= \io v\dix
   +  \int_{\Gamma} v_\Gamma \dS,
\end{equation}
where $v$ represents a generic
function in $L^1(\Omega)\times L^1(\Gamma)$.
With some abuse of notation, we will
also write
\begin{equation}\label{defimm}
  m(v):= \frac1{|\Omega| +  |\Gamma|}
   \ibaro v \dm,
\end{equation}
i.e., the ``mean value'' of $v$ w.r.t.~the measure
$\dm$. Here $|\Gamma|$ represents the surface measure
of $\Gamma$.

With these functional spaces at our disposal, we
can now state our hypotheses on the nonlinear terms. 
%

For convenience, we split $-s_{0,b}'$ (respectively, $-s_{0,\Gamma}'$)
into a sum of a (dominating)
monotone part $f$ (respectively,
$f_\Gamma$) and a linear perturbation.
More precisely, we assume that
\begin{equation}\label{potential_split}
-s_{0,b}'=f(r) -\delta r\quad \perogni r\in\dom f, \quad
 -s_{0,\Gamma}'(r)=f_\Gamma(r) -\delta r \quad \perogni r\in\dom f_\Gamma
\end{equation}
and for some $\delta\geq 0$, with
\begin{equation}\label{compatibility1}
  f \in C^0(\dom f,\RR), \quad
   f_\Gamma \in C^0(\dom f_\Gamma,\RR)
   \quext{monotone,}\quad
   f(0)=f_\Gamma(0)=0,
\end{equation}
where $\dom f$ and $\dom f_\Gamma$ (i.e.,
the {\sl domains}\/ of $f$ and $f_\Gamma$) are
{\sl open}\/ intervals of $\RR$ containing
$0$. We will say that, for instance, $-s_{0,b}$ is 
a ``singular'' potential if its domain 
does not coincide with the whole real 
line (and we will say that it is a ``regular''
potential otherwise). In both cases, we will assume that
\begin{equation}\label{compatibility2}
  \lim_{r\to \de \dom f} (f(r)-\delta r)\sign r
   =\lim_{r\to \de \dom f_\Gamma} (f_\Gamma(r) - \delta r) \sign r
   = +\infty.
\end{equation}
The key assumption that will allow us to obtain a pointwise
estimate of the nonlinear terms $f(\chi)$ and $f\Ga(\chi)$
is the following {\sl compatibility condition}: we ask that
there exist two constants $c_s>0$ and $C_s\ge 0$ such that
\begin{equation} \label{compatibility}
  \dom(f\Ga) \subseteq \dom(f), \quad
  f(r) f\Ga(r) \ge c_s \vert f(r)\vert^2 - C_s
   \,\,\,\perogni r\in \dom(f\Ga).
\end{equation}
In other words, the boundary nonlinear term $f\Ga(r)$, up to the sign,
has to be larger than the bulk nonlinear term $f(r)$, at least for
$r$ far from $0$.
We also introduce, whenever they make sense, the antiderivatives
\begin{equation}\label{defiFG}
  F(r):=\int_0^r f(s)\,\ds, \quad
   F_\Gamma(r):=\int_0^r f_\Gamma(s)\,\ds.
\end{equation}
%
%
%
%
%
Notice that, in case (for instance) $\dom f$ is bounded,
but $f$ is globally summable on $\dom f$,
$F$ can (and will) be extended by continuity
to $\overline{\dom f}$. This is the case, e.g.,
of the logarithmic potential~\eqref{logpot}.
Moreover, $F$ (and, analogously, $F_\Gamma$) will
be thought to be further extended to the whole real line 
by assuming the value $+\infty$ outside its effective domain.
Then, identifying $f$ and $f_\Gamma$ with {\sl maximal monotone graphs}\/
in~$\RR\times\RR$, we have $f=\de F$ and $f_\Gamma=\de F_\Gamma$,
$\de$ representing the {\sl subdifferential}\/ of convex analysis
(here in $\RR$). We refer to the monographs \cite{At,Ba,Br}
for an extensive presentation of the theory
of maximal monotone operators and of
subdifferentials.

To obtain an estimate of the full $V$-norm of $u$,
we will also need a proper form of Poincar\'e's 
inequality (see, e.g., \cite[Lemma 3.2]{ssz12}):
\bele\label{lemma-log-poincare}
 Assume that $\Omega$ is a bounded open subset of $\mathbb{R}^d$. 
 Suppose $v\in W^{1,1}(\Omega)$ and $v\ge 0$ a.e.~in $\Omega$.
 Then, setting $K:=\int_{\Omega}(\log v)^{+}\dx$,
 the following estimate holds:
 \begin{equation}\label{log-poincare}
   \|v\|_{L^1(\Omega)} \le |\Omega| e ^{C_1 K}
    + \frac{C_2}{|\Omega|}\|\nabla v\|_{L^1(\Omega)},
 \end{equation}
 the constants $C_1$ and $C_2$ depending only on $\Omega$.
\enle
\noindent%
Besides assumptions \eqref{potential_split}-\eqref{compatibility}, 
we shall analyze system \eqref{52}-\eqref{55} under the following 
hypotheses:
\begin{align}\label{hyplambda}
  \lambda_b \in C^{2}(\dom f) \hbox{ with } \lambda_b''\in L^\infty(\dom f),
 \quad \lambda_\Gamma \in C^{2}(\dom f_\Gamma) \hbox{ with } \lambda_\Gamma'' \in L^\infty(\dom f_\Gamma),\\
 (\teta_0, \eta_0) \in \calH ,\,\, \teta_0,\eta_0>0\ \text{a.e.},
 \,\,(\log \teta_0, \log\eta_0)\in L^1(\Omega)\times L^1(\Gamma),
 \label{hypteta0}\\
 \chi_0 \in \calV \hbox{ with } s_{0,b}(\chi_0) \in L^1(\Gamma)
\hbox{ and } s_{0,\Gamma}({\chi_{0}}_{\Gamma})\in L^1(\Gamma),\label{hypchi0}\\
%
 %
\mathbf H \in L^2(0,T;\calH), \quad m(\mathbf H)=0
 \ \text{a.e.~in }(0,T).\label{hypH}
\end{align}
Let us define the {\sl energy functional}\/ as
\begin{align} \label{defEnergy}
  \mathcal{E}[(\teta, \eta),(\chi, \chiga)]
   :=&\int_{\Omega} \Big(\teta -\log \teta + \lambda_b(\chi)+\frac{\vert \nabla\chi\vert^2}{2} 
   - s_{0,b}(\chi)\Big)\,\dix\nonumber\\
  & +  \int_{\Gamma} \Big(\eta - \log \eta +\lambda_\Gamma(\chiga) 
         + \frac{\vert \nabla_\Gamma\chiga\vert^2}{2} - s_{0,\Gamma}(\chiga) \Big)\dS,
\end{align}
whenever it makes sense. We will often just write 
$\mathcal{E}[\teta,\chi]$ for brevity.
Note that we admit $\calE$ to take the value $+\infty$.
However, we can readily check that the above assumptions \eqref{hypteta0}-\eqref{hypchi0} 
make the initial energy finite. Namely, we have
\begin{equation} \label{EE0}
  \mathbb{E}_0:=\calE[(\teta_0, \eta_0),(\chi_0, {\chi_0}_{\Gamma})]
   < \infty.
\end{equation}  
Together with \eqref{hyplambda}-\eqref{hypH}, we will also need 
the energy functional~$\calE$ to be coercive with respect
to~$(\chi,\chiga)$. This is obtained
asking that there exist $c_\calV>0$ and $c\in \mathbb{R}$
such that 
\begin{equation} \label{coercivity}
  \mathcal E[\teta,\chi]\ge c_\calV \| \chi\|^2_{\calV} -c
   \ \ \text{for all pairs }(\teta,\chi).
\end{equation}
Such an assumption is satisfied for instance whenever
\begin{equation} \label{compatibilitylat-conf}
  \lambda_{b}(r) -s_{0,b}(r) \ge c_1 r^2 - c_2 
    \hbox{ for all } r\in \dom f\,\,\,\hbox{ and } \,\,\,
  \lambda_{\Gamma}(r)-s_{0,\Gamma}(r) \ge c_1 r^2 -c_2 
    \hbox{ for all } r\in \dom f_\Gamma,
\end{equation}
and for some $c_1>0$ and $c_2\in \mathbb{R}$. Of course, 
\eqref{coercivity} holds trivially if $-s_{0,b}$,
or $-s_{0,\Gamma}$, or both, are singular potentials,
or also in the case of standard double-well potentials.

It is worth noting that the first of assumptions \eqref{hypteta0} 
on $(\teta_0,\eta_0)$ is somehow artificial, in the sense 
that it is stronger than what would be required in order for the
initial energy to be finite. Actually, one could relax \eqref{hypteta0} 
by taking just $(\teta_0, \eta_0) \in L^1(\Omega)\times L^1(\Gamma)$
or even $(\teta_0, \eta_0)\in  \calV'$, paying the price of having
a weaker (and more delicate to deal with) notion of solution.
The Neumann-Neumann Penrose-Fife model with $L^1$ or $V'$ initial 
temperature has been recently studied in \cite{ssz12}. 
In a forthcoming paper, we intend to analyze also the present model
in a similar regularity setting.
%

%
%

%


\subsection{Weak solutions}
\label{subsec:weak}

We use the functional framework introduced above to 
specify a rigorous concept of weak solution to 
the initial value problem for \eqref{56}-\eqref{57},
named Problem~(P) in what follows.
\bedef
 We say that a quadruplet $(\teta, \eta, u, \chi)$
 is a weak solution to~{\rm Problem~(P)} if
 there hold
 \begin{align} \label{regoteta}
   & (\teta, \eta)\in H^1(0,T;\mathcal V')\cap L^\infty (0,T;\calH),\\
  & (\log\teta, \log\eta) \in L^\infty(0,T;L^1(\Omega)\times L^1(\Gamma)),
      \quad \teta, \eta>0\quad\hbox{almost everywhere},\label{positeta}\\
   & u\in L^2(0,T;\mathcal V),\label{regou}\\
   & \chi \in H^1(0,T; \mathcal H)\cap L^\infty(0,T;\mathcal V)\cap L^2(0,T;\mathcal W)\label{regochi},\\
  \label{regofchi}
   & f(\chi) \in \LDH, \quad f\Ga(\chi\Ga) \in L^2(0,T;H\Ga),
 \end{align}
 together with, a.e.~in $(0,T)$, the equations
 \begin{equation}\label{heat}
   {\partial \over {\partial t}}\begin{pmatrix}
    \teta \\ \eta \end{pmatrix}
    = - \mathcal A \begin{pmatrix}
    u \\ u_\Gamma  \end{pmatrix}
   -\begin{pmatrix}\lambda_b(\chi )_t\\ \lambda_\Gamma (\chiga )_t
   \end{pmatrix}
   + \begin{pmatrix}
   h_b \\
   h_\Gamma
  \end{pmatrix}
  \hbox{ in } \calV',
 \end{equation}
 \begin{equation}\label{identif}
   u = -\frac{1}{\teta}\,\,\hbox{ a.e. in }\Omega, \qquad
   \eta = -\frac{1}{u\Ga}\,\,\,\hbox{ a.e. on } \Gamma,
 \end{equation} 
 \begin{equation}\label{phase}
   {\partial \over {\partial t}}\begin{pmatrix}
   \chi  \\ \chiga \end{pmatrix}
   + \mathcal A \begin{pmatrix}
   \chi  \\ \chiga \end{pmatrix}
   - \begin{pmatrix} s'_{0,b}(\chi )\\ s'_{0,\Gamma }(\chiga ) \end{pmatrix}
   = \begin{pmatrix} {{\lambda '_b(\chi )}u }\\ {{\lambda '_\Gamma (\chiga )} u_\Gamma }
   \end{pmatrix} \hbox{ in } \calH,
 \end{equation}
 and the initial conditions
 \begin{equation} \label{CI}
   (\teta, \eta)|_{t=0}=(\teta_0,\eta_0), \quad (\chi, \chiga)|_{t=0}=(\chi_0,{\chi_0}_{\Gamma}) 
    \quad\hbox{a.e.~in } \Omega\ \hbox{and on }\Gamma.
 \end{equation}
 %
 %
 %
\eddef
\noindent%
Sometimes, for brevity, we shall indicate a solution just as a pair
$(\teta,\chi)$ rather than as a quadruplet $(\teta, \eta, u, \chi)$.
\beos\label{tetaeta}
 It is worth giving some explanation on the boundary 
 behavior of $\teta$. Since $u=-1/\teta\in L^2(0,T;\calV)$ 
 by \eqref{regou}, it
 turns out that $u$ has a trace $u\Ga$ on $\Gamma$, 
 which belongs to $V_\Gamma$
 for almost every value of the time variable thanks to
 the definition of $\calV$. On the other hand,  
 we cannot simply write $\eta=\teta\Ga$ since
 the trace of $\teta$ needs not exist. We have, instead,
 to intend $\eta$ as (minus) the reciprocal of the 
 trace of $u$, as specified by the second \eqref{identif}.
 When we consider smoother solutions (for instance, in 
 the a priori estimates, or in the case when 
 $\teta_0$ is more summable, cf.~\eqref{hyptetaLp} below),
 this problem does not occur since the higher regularity of
 $\teta$ permits to give sense to its trace. 
 Regarding the phase variable, the situation is simpler; 
 indeed, by \eqref{regochi}, $\chi$ is always smooth enough 
 to have a trace $\chi\Ga$.
\eddos
\noindent%
%


\subsection{Existence and uniqueness results}
\label{subsec:main}

In this part we state our main results regarding well posedness
of Problem~(P) and regularization properties of weak solutions.
In what follows we will denote by $c$ a positive constant,
which may vary from line to line (or even in the same formula),
depending only on the data of the problem.
Specific dependences will be indicated when needed. 
Moreover, we will denote by $Q$ a nonnegative-valued,
continuous and monotone increasing function of its arguments.
Our main result can be stated as follows:
\bete\label{teor:existenceL2}
 Under assumptions\/ \eqref{potential_split}-\eqref{compatibility},
 \eqref{hyplambda}-\eqref{hypH}, and\/ \eqref{coercivity},
 there exists a unique weak solution to\/ {\rm Problem (P)}.
 Moreover, for all $t>0$ there holds the\/ {\rm energy identity}
 \begin{equation} \label{energy_id}
   \calE[(\teta(t),\eta(t)),(\chi(t),\chi\Ga(t))] 
   + \int_{0}^t \big(\|\nabla u(s)\|^2_{\calH} + \| \chi_t(s)\|_{\calH}^2\big) \ds
    = \mathbb{E}_0 + \itt (\mathbf H(s), u(s))_{\calH} \ds.
 \end{equation}
 If, in addition to~\eqref{hypH},
 the forcing function satisfies
 \begin{equation} \label{hypH+new}
    \mathbf H=(h_b,h_\Gamma)
    \in H^1(0,T;\calV') \cap 
     L^2(0,T,L^{3+\epsilon}(\Omega)\times L^{3+\epsilon}(\Gamma)),
 \end{equation}
 %
 %
 %
 then we also have, for any $\tau\in(0,T)$, the following regularization 
 properties:
 \begin{align}\label{asym-u}
   & \|  u \|_{L^\infty(\tau,T;V)} + \| u \|_{L^\infty((\tau,T)\times \Omega)}
    \le Q( \mathbb{E}_0,\tau^{-1} ),\\
  \label{asym-chi}
   & \| \chi \|_{L^\infty(\tau,T;H^2(\Omega))} + \| f(\chi) \|_{L^\infty(\tau,T;H)}
    \le Q( \mathbb{E}_0,\tau^{-1} ),\\
  & \|  u_\Gamma \|_{L^\infty(\tau,T;V_\Gamma)} + \| u_\Gamma \|_{L^\infty((\tau,T)\times \Gamma)}
    \le Q( \mathbb{E}_0,\tau^{-1} ),  \label{asym-uga}\\
   & \| \chiga \|_{L^\infty(\tau,T;H^2(\Gamma))} + \| f\Ga(\chiga) \|_{L^\infty(\tau,T;H_\Gamma)}
    \le Q( \mathbb{E}_0, \tau^{-1} )  \label{asym-chiga}.
 \end{align}
\ente
\noindent %
The proof of the existence part is based on an approximation-a priori 
bounds-passage to the limit procedure. Additional estimates 
lead to \eqref{asym-u}-\eqref{asym-chiga}.
In particular the $L^\infty$-bounds in \eqref{asym-u} and in 
\eqref{asym-uga} are obtained by adapting a Moser iteration scheme 
with regularization devised in \cite{ssz12}. Note that
the existence part generalizes to the dynamic boundary conditions case the
result of Ito and Kenmochi \cite{Ito-Ken-Kubo03}.

In the case of {\sl singular potentials}, we can also prove that,
at least for strictly positive times, $\chi$ is uniformly separated
from the singularities of the potentials. In order to avoid 
unnecessary technical complications, we shall state this property
under the additional assumption that
\begin{equation} \label{specsing}
  \dom f = \dom f\Ga = (-1,1),~~
   |f\Ga(r)| \ge  \kappa_s |f(r)| - C_s \,\,\, \perogni r\in (-1,1)   
\end{equation}
and for some $\kappa_s\in(0,1]$, $C_s\ge 0$.
Namely, we assume the potentials to be normalized so that the pure
states are represented by the values $\pm1$, both on the bulk
and on the boundary. Moreover, in view of \eqref{compatibility},
we require $|f\Ga|$ to be larger than (a positive constant times)
$|f|$, at least in proximity of $\pm1$. 
\beco\label{cor:sepsing}
 Let the assumptions of\/ {\rm Theorem~\ref{teor:existenceL2}} hold,
 together with\/ \eqref{specsing}. Then, for any $\tau>0$, there 
 exists $\omega\in(0,1)$,
 depending on $\tau$ but independent of $T$, such that 
 \begin{equation} \label{separsing}
   | \chi(t,x) | \le 1 - \omega \ \text{almost everywhere
    in $(\tau,T)\times \Omega$ and in $(\tau,T)\times \Gamma$}.
 \end{equation}
\enco
%
%
 %
 %
%
%
%
\noindent%
Next, we discuss the asymptotic regularization of the temperature field.
Actually, we have already noted (cf.~\eqref{asym-u}, \eqref{asym-uga})
that condition \eqref{hypH+new} on the forcing function 
${\bf H}$ is sufficient in order for $\teta$ to become
uniformly separated from zero for strictly positive times.
The following result (which generalizes \cite[Thm.~2.7]{ssz12},
where Neumann conditions are considered) states that 
$\teta$ is bounded from above, at least for strictly positive
times, provided that \eqref{hypH+new} holds and the initial temperature
enjoys some additional summability property.
%
%
%
%
%
\bepr\label{prop:moser_teta}
 Let the assumptions of\/ {\rm Theorem~\ref{teor:existenceL2}} hold.
 Let also assume that
 \begin{align} \label{hyptetaLp}
   & (\teta_0,\eta_0) \in L^{3+\epsilon}(\Omega)\times L^{3+\epsilon}(\Gamma), 
     \,\,\,\hbox{ for some }\,\,\epsilon>0,\\
  \label{chizbdd}
   & \text{either \eqref{specsing} holds, or }\,\chi_0\in L^\infty(\barO,\dm).
 \end{align}
 Then, any weak solution to\/ {\rm Problem~(P)} satisfies, 
 for any $\tau\in(0,T)$,
 \begin{equation} \label{asym-teta}
  \|\teta \|_{L^\infty((\tau,T)\times \Omega)} 
  + \|\eta \|_{L^\infty((\tau,T)\times \Gamma)}
   \le Q\big(\mathbb{E}_0,\tau^{-1}, \|\teta_0 \|_{L^{3+\epsilon}(\Gamma)}, 
     \|\eta_0 \|_{L^{3+\epsilon}(\Gamma)},
     \|\chi_0 \|_{L^\infty(\barO,\dm)}\big).
 \end{equation}
\empr
\noindent%
We point out that, in the case of singular potentials, the uniform
boundedness of $\chi_0$ required by \eqref{asym-teta}
is a direct consequence of~\eqref{hypchi0}. 

As a byproduct of Prop.~\ref{prop:moser_teta}, we can also prove 
additional regularity of the time derivative of $\teta$:
\beco\label{cor:strong_heat}
 Under the assumptions of\/ {\rm Proposition \ref{prop:moser_teta}}, 
 we have 
 \begin{equation} \label{tetaH1}
   \teta\in H^1(\tau,T;\calH) \,\,\,\, \text{for any}\ \tau>0.
 \end{equation}  
 As a consequence, \eqref{heat} can be decoupled
 and interpreted in the strong form\/ \eqref{57} as a relation
 in $\calH$.
\enco


\section{Proofs of the main results}
\label{sec-proof-exiuni}


\subsection{Proof of Theorem~\ref{teor:existenceL2}: a priori estimates}
\label{subsec:L2}

As a first step, we detail the main estimates constituting
the core of the existence proof. In order to simplify 
the exposition, we limit ourselves to perform 
{\sl formal}\/ a priori bounds on the solutions of 
Problem~(P). In the next section we will see that
these bounds imply weak sequential stability. 
It is clear that, in a formal
proof, the estimates should be performed in the framework
of a proper approximation scheme (e.g., a Faedo-Galerkin
approximation or a time discretization), possibly combined
with some regularization of the data. However, this kind of 
procedure has been already described in full detail in several 
papers related to similar models (see, e.g., \cite{GMS1} and
\cite{gms}) and, actually, the arguments given 
in these papers could be easily adapted to our case. 

That said, we start by presenting the estimates. In all what
follows, we shall assume to have sufficient regularity to justify
all the computations. In particular, we ask $\teta$
to be smooth enough to have a trace, so that 
$\eta=\tetaga$ (and correspondingly $1/\eta=1/\tetaga$).

\smallskip

\noindent%
{\bf Energy estimate.}~~%
Test \eqref{heat} by $\begin{pmatrix}
1 -1/\teta \\ 1-1/\eta \end{pmatrix}$ and \eqref{phase} by
$\begin{pmatrix} \chi_t \\ {\chi_{\Gamma,t}} \end{pmatrix}$.
Noting that two terms cancel out and using that ${\bf H}$ 
has zero mean value, we obtain
\begin{equation} \label{energy_est}
  \ddt \calE[\teta,\chi] + \int_{\Omega}\vert\nabla(-1/\teta)\vert^2\,\dix 
    + \int_{\Gamma}\vert\nabla_\Gamma (-1/\eta)\vert^2\dS 
    + \int_{\Omega}\vert \chi_t\vert^2 \,\dix 
    + \int_\Gamma \vert \chi_{\Gamma,t} \vert^2 \,\dS
  = - ({\bf H}, 1/\teta)_{\calH}.
\end{equation}
Using that $m({\bf H}) = 0$ 
(cf.~\eqref{defimm}), we can write
$$
  ({\bf H}, 1/\teta)_{\calH}
   = \big({\bf H}, 1/\teta-m_{\Omega}\big(1/\teta\big)\big)_{\calH}
   = \int_\Omega h_b \Big(-m_\Omega\Big(\frac{1}{\teta}\Big) 
    + \frac{1}{\teta}\Big) \, \dix
    + \int_\Gamma h_\Gamma \Big(-m_\Omega\Big(\frac{1}{\teta}\Big) 
       + \frac{1}{\eta}\Big)\, \dS. 
$$
Now, the integral over $\Omega$ is easily estimated, using the 
Poincar\'e-Wirtinger inequality, as 
$$  
  \int_\Omega h_b \Big(\frac{1}{\teta} - m_\Omega\Big(\frac{1}{\teta}\Big) \Big)\,\dix 
   \le c_\epsi\|h_b\|^2 
    + \epsi\|\nabla(1/\teta) \|^2 \,\,\,\perogni \epsi>0.
$$
On the other hand, we treat the integral over $\Gamma$ in this way:
\begin{align} \nonumber
  & \disp \int_\Gamma h_\Gamma 
        \Big( \frac{1}{\eta} -m_\Omega\Big(\frac{1}{\teta}\Big) \Big) \,\dS
   \le \|h_\Gamma\|_{\Gamma} \| 1/\eta -m_\Omega(1/\teta) \|_{\Gamma} \\
 \nonumber  
  & \le c \|h_\Gamma \|_{\Gamma} \|1/\teta -m_\Omega(1/\teta)\|_{V} 
  \disp \le c \|h_\Gamma \|_{\Gamma} \big( \|1/\teta -m_\Omega(1/\teta)\|
      + \| \nabla(1/\teta)\| \big) \\
 \label{menomean}
  &  \le c\| h_\Gamma\|_{\Gamma} \| \nabla (1/\teta)\|
   \le c_\epsi \| h_\Gamma\|_{\Gamma}^2 
    + \epsi \| \nabla (1/\teta)\|^2 \,\,\, \perogni \epsi>0,  
\end{align}
where in the second and in the fourth inequalities we have used, 
respectively, the trace theorem and the Poincar\'e-Wirtinger inequality. 
Hence, taking $\epsi$ small enough and integrating \eqref{energy_est}, we obtain
\begin{align}\label{energy_est2}
  \calE(t) + \frac{1}2\int_{0}^t
   \big(\|\nabla(1/\teta)\|^2 
    + \|\nabla\Ga (1/\eta)\|_{\Gamma}^2
    + 2 \| \chi_t\|^2_{\calH} \big)\, \ds
  \le \mathbb{E}_0 + c\|{\bf H}\|^2_{L^2(0,T;\calH)}.
\end{align}
To control the full $V$-norm of $1/\teta$, we have to provide
a bound of its mean value (actually, only the gradient is estimated in
\eqref{energy_est2}). To this aim, we use Lemma~\ref{lemma-log-poincare}
with $v=-1/\teta$, obtaining
\begin{equation}\label{poinc-u}
 \disp \| 1/\teta \|_{L^1(\Omega)} \le \disp |\Omega|e^{c_1 \io(\log\teta)^-}
   + \frac{c_2}{|\Omega|}\|\nabla (1/\teta)\|_{L^1(\Omega)},
%
%
\end{equation}
and the first term on the \rhs\ is uniformly bounded thanks to 
\eqref{energy_est2} and to the expression \eqref{defEnergy}
of the energy functional.
Moreover, by the trace theorem, we get 
\begin{equation}\label{poinc-uga}
   \| 1/ \eta \|_{V_\Gamma} 
    \le c \big( \| 1/ \teta \|_{V} 
          + \| \nabla\Ga (1/ \eta) \|\Ga \big).
\end{equation}
Hence, collecting the above computations, \eqref{energy_est} 
gives the a priori bound
\begin{align} \label{energy_est3}
  \calE[\teta,\chi](t) 
   + \int_{0}^{T} \big( \|(1/\teta)(s) \|^2_{\calV} 
   + \|\chi_t(s) \|^2_{\calH} \big) \, \ds
  \le Q(\mathbb{E}_0,T,\| {\bf H}\|_{L^2(0,T;\calH)}^2)\,\,\,\perogni t\le T.
\end{align}

\smallskip

\noindent%
{\bf Estimate of the nonlinear terms.}~~%
Estimate \eqref{energy_est3} gives that $u=-1/\teta\in L^2(0,T;\calV)$
and $\chi\in L^\infty(0,T;\calV)$, implying, thanks to \eqref{log-poincare},
that the \rhs\ of the phase equation \eqref{phase}
belongs to $L^2(0,T;\calH)$. Now, let us test \eqref{phase} (both on
the bulk and on the boundary) by $f(\chi)$. Using the monotonicity of 
$f$ and assumption \eqref{hyplambda} we then get
\begin{align}\nonumber
  & \ddt \ibaro F(\chi) \,\dm
   + \| f(\chi) \|^2 
   + \iga f(\chi\Ga) \Big( f\Ga(\chi\Ga) 
       - \delta \chi\Ga + \frac{\lambda\Ga'(\chi\Ga)}{\teta\Ga} \Big)\,\dS\\
  \label{giulio11}  
  & \mbox{}~~~~~~~~
   \le \io f(\chi) \Big( \delta \chi - \frac{\lambda_b'(\chi)}{\teta} \Big)\,\dx
   \le \frac12 \| f(\chi) \|^2 
    + c \| \chi \|^2 + c \big(1 + \| \chi \|_{L^4(\Omega)}^2 \big) \| u \|_{L^4(\Omega)}^2.
\end{align}   
The key point is represented by the control of the last term
on the \lhs, and here the compatibility condition \eqref{compatibility}
comes into play. Indeed, using \eqref{compatibility}, 
we readily arrive at
\begin{equation}\label{giulio12}
  \iga f(\chi\Ga) \Big( f\Ga(\chi\Ga) - \delta \chi\Ga + \frac{\lambda\Ga'(\chi\Ga)}{\teta\Ga} \Big)
   \ge \frac{c_s}2 \| f(\chi\Ga) \|\Ga^2 - c - c \| \chi\Ga \|\Ga^2 
      - c \big( 1 + \| \chi\Ga \|_{L^4(\Gamma)}^2 \big) \| u\Ga \|_{L^4(\Gamma)}^2.
\end{equation}
Hence, integrating \eqref{giulio11} in time and using \eqref{giulio12} and 
\eqref{energy_est3}, we infer
\begin{equation}\label{giulio13}
  \| f(\chi) \|_{L^2(0,T;H)} 
   + \| f(\chi\Ga) \|_{L^2(0,T;H\Ga)} \le Q(\EE_0,T).   
\end{equation}
Here and below, we allow $Q$ to depend additionally on ${\bf H}$.
Moreover, we recall that $\| F(\chi_0) \|_{L^1(\barO,\dm)} \le Q(\EE_0)$ 
thanks to \eqref{hypchi0} and \eqref{compatibility}.

Next, \eqref{energy_est3}, \eqref{giulio13}, and a comparison 
of terms in \eqref{52} (i.e., the bulk component
of \eqref{heat} in the strong formulation -- recall that
we assume the solutions to be smooth at this level), give
\begin{equation}\label{cons_energy2}
   \|\Delta\chi \|_{L^2(0,T;H)}
     \le Q(\mathbb{E}_0,T).
\end{equation}
Consequently, using standard trace and elliptic regularity theorems
(cf., e.g.~\cite[Theorem 2.7.7 and Theorem~3.1.5]{BG}),
it is not difficult to arrive at 
\begin{equation}\label{giulio14}
  \| \dn\chi \|_{L^2(0,T;H\Ga)} \le Q(\mathbb{E}_0,T).
\end{equation}
This allows to test the boundary equation \eqref{53} 
by $f\Ga(\chi\Ga)$. Proceeding as above
(but {\sl without}\/ using \eqref{52}) and 
controlling the term $\dn \chi$ directly by means
of \eqref{giulio14}, we finally arrive at
\begin{equation}\label{giulio15}
  \| \Delta\Ga\chi\Ga \|_{L^2(0,T;H\Ga)} 
   + \| f\Ga(\chi\Ga) \|_{L^2(0,T;H\Ga)} \le Q(\mathbb{E}_0,T).
\end{equation}

\smallskip

\noindent%
Now, let us come to the temperature equation. 
By \eqref{energy_est3} and assumption \eqref{hyplambda},
we have 
\begin{equation}\label{giulio21}
  \| \partial_t \lambda_b(\chi) \|_{L^2(0,T;L^{3/2}(\Omega))}
   \le \| \lambda_b'(\chi) \|_{L^\infty(0,T;L^6(\Omega))} \| \chi_t \|_{\LDH}
   \le Q(\mathbb{E}_0,T),
\end{equation}
and a similar relation on $\Gamma$. 
Actually, when $-s_{0,b}$ and $-s_{0,\Gamma}$
are singular potentials, we automatically have a uniform 
$L^\infty$-bound on $\chi$.
Hence, we obtain more precisely
\begin{equation}\label{giulio22}
  \| \partial_t \lambda_b(\chi) \|_{\LDH}
   \le \| \lambda_b'(\chi) \|_{L^\infty(Q)} \| \chi_t \|_{\LDH}
   \le Q(\mathbb{E}_0,T).
\end{equation}
By a comparison of terms in the (coupled) 
weak formulation \eqref{heat}, we then get
in both cases
\begin{equation}\label{giulio23}
  \| \teta_t \|_{L^2(0,T;\calV')}
   \le Q(\mathbb{E}_0,T).
\end{equation}
Finally, to get the $\calH$-regularity of $\teta$, we test \eqref{heat}
by $\teta$. Then, using the first of \eqref{hypteta0}, \eqref{hypH} and 
\eqref{giulio22}, noting that 
\begin{equation}\label{betti11}
  \bigg| \io \lambda_b'(\chi) \chi_t \teta \bigg|
   \le c \big( 1 + \| \chi \|_{L^\infty(\Omega)} \big)
    \| \chi_t \| \| \teta \|
   \le c \| \chi_t \|^2
   + c \big( 1 + \| \chi \|_{H^2(\Omega)}^2 \big)
     \| \teta \|^2,
\end{equation}
thanks also to \eqref{hyplambda}, observing that a similar 
relation holds on $\Gamma$, and applying Gronwall's lemma,
it is not difficult to arrive at
\begin{equation}\label{L2teta}
  \| \teta \|_{L^\infty(0,T;\calH)}\le Q(\mathbb{E}_0,T).
\end{equation}
Relations \eqref{giulio23} and \eqref{L2teta}
give the desired \eqref{regoteta}.
%
%
%
%
%


\subsection{Proof of Theorem~\ref{teor:existenceL2}: weak sequential stability}

\label{sec:passlim}

In this part, we prove weak sequential stability for
Problem~(P). Namely, we assume to have a sequence $(\teta_n,\chi_n)$ 
of sufficiently smooth solutions having uniformly bounded initial 
energies~$\EE_{0,n}$ and show that, at least up to the
extraction of subsequences, this family converges in a 
suitable sense to a weak solution to the problem. 
Of course, in principle $(\teta_n,\chi_n)$ could also be
taken as the solution to some regularized Problem~(P$_n$);
in this respect, the procedure below can be seen as a simplified 
version of the arguments needed to remove the regularization.
However, as noticed at the beginning of Subsec.~\ref{subsec:L2}, 
we decided, for the sake of simplicity, to skip the details of 
the approximation procedure.

That said, set $u_n :=- 1/ \teta_n $. 
Consequently, its trace $u_{n,\Gamma}$ 
is equal to $-1/\eta_n$. Thanks to estimates \eqref{energy_est3}, \eqref{cons_energy2},
\eqref{giulio15} and \eqref{giulio23}-\eqref{L2teta},
and to the assumed uniform boundedness of the initial energies,
we then obtain
\begin{align} \label{pass_lim1}
  \sup_{n\in\mathbb{N}}& \left\{\| u_n \|_{L^2(0,T;\calV)} 
   + \|\partial_t\chi_n \|_{L^2(0,T;\calH)} 
   + \| \chi_n \|_{L^2(0,T;\calW)} 
    + \|s'_{0,b}(\chi_n) \|_{L^2(0,T;H)}\right. \\
  \nonumber 
  & \left. + \|s'_{0,\Gamma}(\chi_{n,\Gamma} ) \|_{L^2(0,T;H_\Gamma)} 
   + \| \partial_t \teta_n \|_{L^2(0,T;\calV')}
   + \| \teta_n \|_{L^\infty(0,T;\calH)}
        \right\} 
  \le c,
\end{align}
where $c$ is a positive constant independent of $n$. Then, using
standard weak compactness arguments and the Aubin-Lions lemma,
there exist $u, \chi$ with traces, respectively, $u_{\Gamma},\chiga$, and a 
non-relabelled subsequence of $n$ such that
\begin{align}
 &\disp u_n \xrightarrow{n\nearrow \infty} u \hbox{ weakly in } L^2(0,T;\calV),\label{weak_u}\\
 &\disp  \chi_n \xrightarrow{n\nearrow \infty} \chi \hbox{ weakly in } H^1(0,T;\calH)\cap L^2(0,T;\calW)
  ~\hbox{ and strongly in } L^2(0,T; \mathcal V),
  \label{weak_chi}\\
 &\disp s'_{0,b}(\chi_n)\xrightarrow{n\nearrow \infty} \overline{s'_{0,b}} \hbox{ weakly in } L^2(0,T;H),
  \label{weak_sb}\\
 &\disp  s'_{0,\Gamma}(\chi_n)\xrightarrow{n\nearrow \infty} \overline{s'_{0,\Gamma}}
 \hbox{ weakly in }   L^2(0,T;H_\Gamma),\label{weak_sGa}\\
 &\disp \lambda'_{b}(\chi_n) u_n \xrightarrow{n\nearrow \infty} \overline{\lambda'_{b}} \hbox{ weakly in }
    L^2(0,T;H),\label{weak_lambdab}\\
 &\disp \lambda'_{\Gamma}(\chi_{n,\Gamma}) u_{n,\Gamma} \xrightarrow{n\nearrow \infty}
  \overline{\lambda'_{\Gamma}} \hbox{ weakly in } L^2(0,T;H_\Gamma)\label{weak_lambdaGa}.
\end{align}
Then, thanks to the properties of $\lambda$ and $\lambda_\Gamma$, 
it is not difficult to obtain that $\lambda'(\chi_n)\longrightarrow \lambda'(\chi)$ 
strongly in $L^2(0,T;H)$ (similarly for $\lambda'_{\Gamma}$).
Moreover, by the usual monotonicity argument \cite[Prop.~1.1, p.~42]{Ba} 
(recall that $s_{0,b}'$ and $s_{0,\Gamma}'$ verify \eqref{potential_split}),
we can identify, almost everywhere in $(0,T)\times \Omega$ 
(or in $(0,T)\times \Gamma$), $\overline{s'_{0,b}}=s'_{0,b}(\chi)$, 
$\overline{s'_{0,\Gamma}}=s'_{0,\Gamma}(\chiga)$, $\overline{\lambda'_{b}}=\lambda'_{b}(\chi) u$ 
and $\overline{\lambda'_{\Gamma}}=\lambda'_{\Gamma}(\chiga)u_\Gamma$.
Now, we can take the limit in \eqref{heat}. Indeed, by 
\eqref{L2teta} we have
\begin{equation}\label{conf11}
  (\teta_n,\teta_{n,\Gamma}) \xrightarrow{n\nearrow \infty} (\teta,\eta) \hbox{ weak star in } 
   L^\infty(0,T;\calH).
\end{equation}
%
To conclude we have to show that $-1/\teta= u$ 
a.e.~in $(0,T)\times \Omega$ and that $-1/\eta=u_\Gamma$ a.e.~on $(0,T)\times \Gamma$.
Actually, these identifications are consequences of the monotonicity 
of the map $v\mapsto -1/v$ (and of its realization in 
$L^2(0,T;\calH)$). Indeed, by~\eqref{weak_u},
\begin{equation} \label{conf12}
  \disp  u_n \xrightarrow{n\nearrow \infty} u \hbox{ weakly in } L^2(0,T;\calV)
   \hbox{ and consequently in }L^2(0,T;L^6(\Omega)\times L^6(\Gamma)),
\end{equation}
Moreover, noting that, thanks to \eqref{giulio23} and \eqref{L2teta},
\begin{equation}\label{giulio31}
  \teta_n\xrightarrow{n\nearrow \infty} \teta 
   \hbox{ strongly in } L^2(0,T;\calV'),
\end{equation}
we get the following limsup-inequality:
\begin{equation}\label{giulio32}
  \limsup_{n\nearrow \infty} \int_{0}^T \ibaro \teta_n u_n \, \dm \, \dt
   = \limsup_{n\nearrow \infty} \int_{0}^T \duav{\teta_n,u_n}_{\calV} \, \dt
   = \int_{0}^T \duav{\teta,u}_{\calV} \, \dt
   = \int_{0}^T \ibaro \teta u \, \dm \, \dt.
\end{equation}
Combining \eqref{conf11}-\eqref{conf12} with \eqref{giulio32}
and using as above \cite[Prop.~1.1, p.~42]{Ba},
we obtain $-1/\teta= u$ a.e.~in $(0,T)\times \Omega$
and $-1/\eta=u_\Gamma$ a.e.~on $(0,T)\times \Gamma$.
This concludes the proof of weak sequential
stability. In particular, this argument could be easily adapted to 
provide a rigorous existence proof. 

Finally, the energy identity \eqref{energy_id} can be proved simply by
multiplying \eqref{heat} by $1-1/\teta$ and \eqref{phase} by $\chi_t$
and performing standard integrations by parts.
Indeed, it is easy to check that, also in the limit, 
$1-1/\teta$ and $\chi_t$ are sufficiently smooth 
to be used as test functions.


\subsection{Proof of Theorem~\ref{teor:existenceL2}: uniqueness and regularity}
\label{subsec:rego}

{\bf Uniqueness.}~~%
Uniqueness can be easily obtained as follows: take two weak solutions 
$(\teta_1, \chi_1)$ and $(\teta_2, \chi_2)$ and set 
$(\teta,\chi):=(\teta_1-\teta_2,\chi_1-\chi_2)$. Then, write 
\eqref{heat} firstly for $(\teta_1,\chi_1)$, then for $(\teta_2, \chi_2)$, take
the difference, integrate it over $(0, t)$, $t\in (0, T]$, and test the 
result by $u_1-u_2$, where $u_i=-1/\teta_i$. Correspondingly,
take the difference of \eqref{phase} and 
test it by $\chi$. Using the monotonicity of $f$ and $f\Ga$ 
and proceeding as in \cite[Sec.~3]{RSMed} in order to 
estimate the nonlinearities coming from the quadratic terms 
involving $\lambda_b$ and $\lambda_\Gamma$, it is not difficult
to achieve the desired result. The details, very similar to 
the proof given in \cite[Sec.~3]{RSMed}, are left to the reader.

\smallskip

\noindent%
{\bf Parabolic regularization of solutions.}~~%
Now we come to the proof of~\eqref{asym-u}-\eqref{asym-chiga}. 
As before, we proceed by formal estimates. As noted above, 
also these estimates should be performed within
some approximation scheme, whose details are omitted for
simplicity.

\smallskip

\noindent%
{\bf Second estimate.}~~%
We test \eqref{heat} by $t u_t=t \teta_t/\teta^2$
and add the result to  the time derivative of \eqref{phase}
multiplied by $t \chi_t$.
%
Recalling \eqref{potential_split} and using the monotonicity of $f$ and 
of $f_\Gamma$,
we obtain
\begin{align}\nonumber
  & \disp\ddt \Big( \frac{t}2 \| \nabla u \|^2 + \frac{t}2\|\nabla_\Gamma u_\Gamma\|^2_\Gamma
   + \frac{t}2 \| \chi_t \|^2_\calH \Big)
   + t\io \teta_t^2 u^2 + t\iga\eta_t^{2} u_\Gamma^{2}
   + t \| \nabla \chi_t \|^2 + t\|\nabla_\Gamma \chi_{\Gamma,t} \|^2_\Gamma\\
  \nonumber
  & \disp \mbox{}~~~~~
   - t\io h_b u_t - t\iga h_\Gamma u_{\Gamma,t}
   \le \frac12\big( \| \nabla u \|^2 +\| \nabla_\Gamma u_\Gamma \|^2_\Gamma\big)
   + \Big(\frac12 +\delta\Big) \| \chi_t\|^2_{\calH} \\
 \label{conto31}  
  & \mbox{} ~~~~~~~~~~
   + t\iga \lambda''_{\Gamma}(\chiga) \chi_{\Gamma,t}^2 u_\Gamma
   + t \io\lambda''_b(\chi)\chi_t^{2} u.
\end{align}
Now, let us note that, by \eqref{hyplambda}, 
we have for any $\epsi>0$
\begin{equation} \label{conto31bis}
  t \io\lambda''(\chi)\chi_t^{2} u \le t\Big(\epsi \| \chi_t\|^2_{V} 
   + c_\epsi \|\chi_t \|^2 \| u \|^2_{V} \Big),
   \quad t \iga\lambda_\lambda''(\chiga)\chi_{\Gamma,t}^{2} u_\Gamma 
        \le t\Big(\epsi \| \chi_{\Gamma,t} \|^2_{V_\Gamma} 
     + c_\epsi \|\chi_t \|^2_\Gamma \| u_\Gamma \|^2_{V_\Gamma} \Big).
\end{equation}
Moreover,
\begin{equation} \label{co31c}
  - t \io h_b u_t - t \iga h\Ga \ugat
  = - \ddt \Big( \io t u h_b + \iga t u\Ga h\Ga \Big)
    + \Big( \io u h_b + \iga u\Ga h\Ga \Big)
    + \Big( \io t u \de_t h_b + \iga t u\Ga \de_t h\Ga \Big).
\end{equation}
Now, thanks to the fact that ${\bf H}$ has zero mean value
(cf.~\eqref{hypH}), we can modify all the integrands on the 
\rhs\ by subtracting to $u$ and $u\Ga$ the quantity $m\OO(u)$,
as done in \eqref{menomean}. For instance, the first pair of
integrals gives
\begin{equation} \label{co31d}
  - \ddt \Big( \io t u h_b + \iga t u\Ga h\Ga \Big)
    = - \ddt \Big( \io t ( u - m\OO(u) ) h_b 
         + \iga t ( u\Ga - m\OO(u) ) h\Ga \Big).
\end{equation}
The other two pairs of integrals on the \rhs\ of \eqref{co31c} are 
managed similarly and then estimated directly in this way:
\begin{equation}\label{co31e}   
  \Big( \io u h_b + \iga u\Ga h\Ga \Big)
    + \Big( \io t u \de_t h_b + \iga t u\Ga \de_t h\Ga \Big)
   \le c \| u \|_{\calV} \| {\bf H} \|_{\calV'}
    + c t \| u \|_{\calV} \| \partial_t {\bf H} \|_{\calV'}.
\end{equation}
%
%
%
%
%
As a consequence, taking $\epsi$ small enough in \eqref{conto31bis},
integrating \eqref{conto31} between $0$ and a generic 
$t\in (0,T)$, recalling \eqref{hypH+new} 
and \eqref{energy_est3}, and using Gronwall's lemma, it is not difficult 
to arrive at
\begin{align} \nonumber
  & \disp \frac{t}2\| \nabla u (t) \|^2 
  + \frac{t}2\|\nabla_\Gamma u_\Gamma(t)\|^2_\Gamma
  + \frac{t}2 \| \chi_t (t) \|^2_\calH 
  - \io t \big( u(t) - m\OO(u(t)) \big) h_b(t) \\
 \nonumber  
  & \mbox{}~~~~~
  - \iga t \big( u\Ga(t) - m\OO(u(t)) \big) h\Ga(t)  
   + \int_{0}^t s\io \teta_t^2 u^2 
   + \int_{0}^t s\iga\eta_t^{2} u_\Gamma^{2} \\
 \label{local_est} \disp
  & \mbox{}~~~~~ 
   + \frac12 \int_{0}^t s \| \nabla \chi_t \|^2 
   + \frac12 \int_0^{t} s\|\nabla_\Gamma \chi_{\Gamma,t} \|^2_\Gamma
  \le Q(\mathbb{E}_0),
\end{align}
where, as before, the expression of $Q$ may depend on the 
source ${\bf H}$ (and in particular on the additional 
regularity assumptions \eqref{hypH+new}). Note 
now that
\begin{equation}\label{st-inte}
  - \io t \big( u(t) - m\OO(u(t)) \big) h_b(t) 
  - \iga t \big( u\Ga(t) - m\OO(u(t)) \big) h\Ga(t) 
   \ge - \disp \frac{t}4 \| \nabla u (t) \|^2 
   - \frac{t}4\|\nabla_\Gamma u_\Gamma(t)\|^2_\Gamma
   - c, 
\end{equation}
where $c=c({\bf H})$ is independent of $t$. Hence,
using again the generalized Poincar\'e inequality
of Lemma~\ref{lemma-log-poincare}, we finally obtain
\begin{equation}\label{global_est2}
  \| u \|^2_{L^\infty(\tau,T,\calV)} 
   + \| \chi_t \|^2_{L^\infty(\tau,T;\calH)} + \|\chi_t \|^2_{L^2(\tau,T;\calV)}
  \le Q(\mathbb{E}_0,\tau^{-1},T),
\end{equation}
which gives the first of \eqref{asym-u} and of \eqref{asym-uga}. 
Moreover, testing again \eqref{phase}
by $f(\chi)$ we can write (compare with \eqref{giulio11})
\begin{align}\nonumber
  & \| f(\chi) \|^2 
   + \iga f(\chi\Ga) \Big( f\Ga(\chi\Ga) 
       - \delta \chi\Ga + \frac{\lambda\Ga'(\chi\Ga)}{\teta\Ga} \Big) \\
 \nonumber
  & \mbox{}~~~~~
   \le \io f(\chi) \Big( - \chi_t + \delta \chi - \frac{\lambda_b'(\chi)}{\teta} \Big)
    - \iga f(\chi\Ga) \chi_{\Gamma,t} \\
 \label{giulio11x}  
  & \mbox{}~~~~~
   \le \frac12 \| f(\chi) \|^2 
   + \frac{c_s}4 \| f(\chi\Ga) \|^2\Ga
   + c \big( \| \chi \|^2 + \| \chi_t \|^2 + \| \chi_{\Gamma,t} \|^2\Ga \big)
   + c \big(1 + \| \chi \|_{L^4(\Omega)}^2 \big) \| u \|_{L^4(\Omega)}^2.
\end{align}   
Consequently, recalling \eqref{giulio12}, taking the essential
supremum with respect to time, and using \eqref{hyplambda}, 
\eqref{global_est2}, \eqref{energy_est3}, and the 
compatibility condition \eqref{compatibility},
it is not difficult to infer
%
%
%
\begin{equation} \label{pot-LinfH}
  \|-\Delta \chi\|_{L^\infty(\tau,T;H)} 
   + \|s'_{0,b}(\chi)\|_{L^\infty(\tau,T;H)}
  \le Q(\mathbb{E}_0,\tau^{-1},T),
   \hbox{ for any } \tau\in (0,1).
\end{equation}
Then, using again \cite[Theorem 2.7.7 and Theorem 3.1.5]{BG},
we get $\dn \chi \in L^\infty(\tau,T;H_\Gamma)$ 
for any $\tau\in (0,1)$. Hence, following the lines 
of the argument in Section~\ref{subsec:L2},
we arrive at (compare with \eqref{giulio14}-\eqref{giulio15}) 
\begin{equation} \label{pot-LinfHga}
  \|-\Delta_\Gamma \chiga\|_{L^\infty(\tau,T;H_\Gamma)} 
    + \|s'_{0,\Gamma}(\chiga)\|_{L^\infty(\tau,T;H_\Gamma)}
  \le Q(\mathbb{E}_0,\tau^{-1},T),
   \hbox{ for any } \tau\in (0,1).
\end{equation}
Relations~\eqref{pot-LinfH}, \eqref{pot-LinfHga} give
\eqref{asym-chi}, \eqref{asym-chiga}, respectively. Hence,
to conclude the proof of Theorem~\ref{teor:existenceL2},
it remains to show the second of \eqref{asym-u} and of
\eqref{asym-uga}. This essentially relies on the following
result, which is a variant of \cite[Lemma 3.3]{ssz12}:
\bele \label{lemma:moser-u}
 Let $(\teta, \eta)$ be a smooth solution to the problem
\begin{align}\label{mos_u}
  & \partial_t\teta - \Delta u = g_b,
   \quad u = -\frac{1}{\teta},
   \quext{in } \Omega,\\
   & \partial_t \eta -\Delta_\Gamma u_\Gamma = g_\Gamma - \dn u,
    \quad \eta = -\frac{1}{u_\Gamma}\,
  \quext{on }\Gamma, \label{mos_uga}
 \end{align}
 over the generic time interval $(S,T)$, where we additionally
 assume that
 \begin{equation}\label{inpiu}
   \| u \|_{L^3(S,T;L^{3/2}(\Omega))} \le M, \qquad
    \| g_b \|_{L^2(S,T;L^{3+\epsilon}(\Omega))} 
     + \| g_\Gamma \|_{L^2(S,T;L^{3+\epsilon}(\Gamma))} \le G,
 \end{equation}
 for some (given) constants $M>0$, $G>0$ and some
 $\epsilon>0$. Moreover, let us assume that
 \begin{equation}\label{uSL1}
   u(S)\in L^1(\Omega),\quad u_\Gamma(S)\in L^1(\Gamma), \qquad
    \| u(S) \|_{L^1(\Omega)} + \| u\Ga(S) \|_{L^1(\Gamma)} \le U,
 \end{equation}
 for some $U>0$. Then, we have
 \begin{equation}\label{st-abstr-u}
   \| u \|_{L^\infty((S+\tau,T)\times \Omega)} 
     + \| u_\Gamma \|_{L^\infty((S+\tau,T)\times \Gamma)} \
    \le Q\big( G, M, U, \tau^{-1}\big) \,\,\,\perogni \tau \in (0,T-S).
 \end{equation}
\enle
\noindent%
%
The proof is mainly based on a Moser iteration scheme with 
regularization and closely follows the lines of \cite[Lemma 3.3]{ssz12},
where the case of Neumann boundary conditions is treated. The easy adaptation
is left to the reader.

At this point we are in a position to conclude the 
proof of Theorem~\ref{teor:existenceL2}. Actually, 
assumption \eqref{hypH+new} and estimate \eqref{global_est2}
allow us to apply Lemma~\ref{lemma:moser-u} with 
$g_b = h_b - (\lambda_b(\chi))_t$ and 
$g_\Gamma = h_\Gamma - (\lambda_\Gamma(\chi\Ga))_t$
on the generic interval $(S,T)$ with $S>0$.
This provides the $L^\infty$-regularization for $u$ and for 
its trace $u_\Gamma$, as desired.~~\dimbox

\smallskip

\noindent%
{\bf Proof of Corollary~\ref{cor:sepsing}.}~~%
The proof is based on the comparison principle for ODEs, similarly to
\cite[Sec.~3.2]{GPS1}.
Given $\tau>0$, thanks to the bounds \eqref{asym-u}-\eqref{asym-chiga} and to 
assumption \eqref{hyplambda}, there exists $M>0$ depending on $\EE_0$ and 
$\tau$, but otherwise independent of time, such that
\begin{equation} \label{giulio42}
  \| \delta\chi + \lambda_b'(\chi) u \|_{L^\infty((\tau,T)\times\Omega)}
  + \| \delta\chi\Ga + \lambda\Ga'(\chi\Ga) u\Ga \|_{L^\infty((\tau,T)\times\Omega)}
  \le M.
\end{equation}
Let us now consider the following (forward in time)
initial value problem
\begin{equation} \label{giulio43}
  \begin{cases}
    y' + \kappa_s f(y) = M + C_s,\\
    y(\tau) = 1,
  \end{cases} 
\end{equation}
the constants $\kappa_s$ and $C_s$ being as in \eqref{specsing}.
Thanks to \eqref{compatibility2} and to
the standard theory of ODE's, there 
exists $\omega=\omega(\tau)\in (0,1)$ such that 
$y(t) \le 1-\omega$ for all $t\ge 2 \tau$. Let us then view
$y$ as a function also of the space variable $x$ and 
subtract the first row of \eqref{giulio43} from \eqref{phase},
both in the bulk and on the boundary. Then, test the result
by $(\chi-y)^+$, $(\cdot)^+$ denoting the positive part. We have
\begin{align} \nonumber
  & \frac12 \ddt \ibaro \big| (\chi-y)^+ \big|^2 \dm
   + \io \big( f(\chi) - \kappa_s f(y) \big) (\chi-y)^+ \dx\\
 \nonumber
  & \mbox{}~~~~~  
   + \iga \big( f\Ga(\chi\Ga) - \kappa_s f(y) + C_s \big) (\chi\Ga-y)^+ \dS\\
 \label{giulio41}  
  & \mbox{}~~~~~  
   = \io \big( \delta\chi + \lambda_b'(\chi) u - M - C_s \big) (\chi-y)^+  \dx 
   + \iga \big( \delta\chi\Ga + \lambda\Ga'(\chi\Ga) u\Ga - M \big) (\chi\Ga-y)^+ \dS
   \le 0.
\end{align}   
Recalling \eqref{specsing} and using in particular that $\kappa_s\in(0,1]$,
exploiting the monotonicity of $f$, and noting that $(\chi-y)^+|_{t=\tau}=0$ 
$\dm$-almost everywhere in $\barO$ thanks to the fact that 
$\chi$ takes almost everywhere
its values in $\dom f=(-1,1)$, the comparison principle yields
that $\chi(t,x)\le y(t) \le 1 - \omega$ for all $t\ge 2\tau$
and for $\dm$-almost all $x\in\barO$. Namely, we have the upper bound
in \eqref{separsing}. The lower bound is proved similarly.~~\dimbox
%
 %
 %
%

\smallskip

\noindent%
{\bf Proof of Prop.~\ref{prop:moser_teta}.}~~%
As before, we get the uniform boundedness of $\teta$ by
means of a Moser iteration scheme. Namely, we will
rely on the following lemma, whose proof can be obtained
by suitably modifying the argument given in
\cite[Lemma~3.5]{ssz12} and is thus left
to the reader.
\bele\label{moser-teta}
 Let $(\teta,\eta)$ be a smooth solution of problem\/
 \eqref{mos_u}-\eqref{mos_uga} over the generic time 
 interval $(S,T)$, where we additionally assume that
 \begin{equation}\label{inpiu-2}
   \| \teta \|_{L^\infty(S,T;L^1(\barO,\dm))} \le M, \qquad
    \| g_b \|_{L^2(S,T;L^{3+\epsilon}(\Omega))} 
     + \| g_\Gamma \|_{L^2(S,T;L^{3+\epsilon}(\Gamma))} \le G,
 \end{equation}
 for some (given) constants $M>0$, $G>0$ and some
 $\epsilon>0$. Moreover, let us assume that 
 \begin{equation}\label{inpiu-2b}
   \teta(S)\in L^{3+\epsilon}(\Omega),~~\eta(S)\in L^{3+\epsilon}(\Gamma),
    \qquad \| \teta(S) \|_{L^{3+\epsilon}(\Omega)}
    + \| \eta(S) \|_{L^{3+\epsilon}(\Omega)} \le \Theta,
 \end{equation}
 for some (given) constant $\Theta>0$.
 Then, for any $\tau \in (0,T-S)$, we have
 \begin{equation}\label{st-abstr-teta}
   \| \teta \|_{L^\infty((S+\tau,T)\times \Omega)} \le 
    Q\big( M, G, \Theta, \tau^{-1} \big).
 \end{equation}
\enle
\noindent%
However, to apply the above lemma over the generic interval $(S,T)$,
$S>0$, some care is needed. Actually, taking
as above $g_b = h_b - (\lambda_b(\chi))_t$ and 
$g_\Gamma = h_\Gamma - (\lambda_\Gamma(\chi))_t$,
\eqref{inpiu-2} is satisfied
thanks to hypothesis \eqref{hypH+new} and
estimate \eqref{global_est2}. However, it is not 
a priori obvious that assumption \eqref{hyptetaLp}
implies \eqref{inpiu-2b} for $S>0$.
To prove this fact, we have
to provide a uniform control of the 
$L^{3+\epsilon}$-norm of $\teta$ over the interval $(0,S)$.

To this aim, we just consider the situation
(more difficult, here) when $-s_{0,b}$ and
$-s_{0,\Gamma}$ are ``regular''
potentials and \eqref{chizbdd} holds.
In this case, we first have to prove that a uniform control
of the $L^p$-norm of $\chi$ holds on small time intervals
(this property is, of course, obvious when we have singular
potentials). To do this, it suffices
to test \eqref{phase} by $|\chi|^{p-2}\chi$, $p$ to be
chosen below, and integrate 
both in $\Omega$ and on $\Gamma$. Using the monotonicity of 
$f$ and $f\Ga$ together with assumption
\eqref{hyplambda}, it is then easy to get, for some $\kappa>0$,
\begin{align} \nonumber
  & \ddt \| \chi \|^{p}_{L^{p}(\barO,\dm)}
   + \kappa \io \big| \nabla ( |\chi|^{\frac{p-2}2} \chi ) \big|^2
   + \kappa \iga \big| \nabla\Ga ( |\chi\Ga|^{\frac{p-2}2} \chi\Ga ) \big|^2\\
 \label{contok1}
  & \mbox{}~~~~~
   \le c \ibaro |\chi|^p \dm
    + c \ibaro \big( 1 + |\chi|^p \big) |u| \,\dm
\end{align}
and the latter term on the \rhs\ can be controlled
as follows (we just deal with the bulk component, the
boundary one behaving similarly):
\begin{align}\nonumber
  \io  |\chi|^p |u| 
   & \le \| u \|_{L^6(\Omega)} 
    \| \chi \|_{L^\frac{6p}5(\Omega)}^p
   \le \| u \|_{L^6(\Omega)} 
    \| \chi \|_{L^p(\Omega)}^{\frac{3p}4}
    \| \chi \|_{L^{3p}(\Omega)}^{\frac{p}4} 
   \le c \| u \|_{L^6(\Omega)} 
    \| \chi \|_{L^p(\Omega)}^{\frac{p}2}
    \| \chi \|_{L^{3p}(\Omega)}^{\frac{p}2} \\
 \label{contok2}   
   & \le c_\sigma  \| u \|_{V}^2
    \| \chi \|_{L^p(\Omega)}^p
    + \sigma \| \chi \|_{L^{3p}(\Omega)}^p
   \le c_\sigma  \big( 1 + \| u \|_{V}^2 \big)
    \| \chi \|_{L^p(\Omega)}^p
    + \sigma \big\| \nabla (|\chi|^{\frac{p-2}2} \chi) \big\|^2.
\end{align}
Then, we first take $\sigma$ small enough so that
the latter term on the \rhs\ is controlled by 
the second term on the \lhs\ of \eqref{contok1}. Subsequently,
we use Gronwall's Lemma in \eqref{contok1} to obtain,
for any $p\in(1,\infty)$,
\begin{equation}\label{contok3}
  \| \chi \|_{L^\infty(0,S;L^p(\barO,\dm))} 
   \le Q\big(\EE_0,\| \chi_0 \|_{L^\infty(\barO,\dm)},S,p\big).
\end{equation}   
Next, we proceed along the lines of 
\cite[Proof of Theorem 2.7]{ssz12}.
Namely, we test \eqref{heat} by $\teta^{2+\epsilon}$,
both on the bulk and on the boundary.
We then get
\begin{align}\nonumber
  \ddt \| \teta \|^{3+\epsilon}_{L^{3+\epsilon}(\barO,\dm)}
   & \le c \io | - \lambda_b'(\chi)\chi_t + h_b | \teta^{2+\epsilon}
    + c \iga | - \lambda\Ga'(\chi\Ga)\chi_{\Gamma,t} + h\Ga | \eta^{2+\epsilon}\\
 \label{contox1}
  & \le c \Big(
   \| - \lambda_b'(\chi)\chi_t + h_b \|_{L^{3+\epsilon}(\Omega)}
    + \| - \lambda\Ga'(\chi\Ga)\chi_{\Gamma,t} + h\Ga \|_{L^{3+\epsilon}(\Gamma)} \Big)
    \| \teta \|^{2+\epsilon}_{L^{3+\epsilon}(\barO,\dm)},
\end{align}   
whence, clearly,
\begin{equation}\label{contox2}
  \ddt \| \teta \|_{L^{3+\epsilon}(\barO,\dm)}
   \le c \Big( \| \lambda_b'(\chi)\chi_t \|_{L^{3+\epsilon}(\Omega)}
     + \| h_b \|_{L^{3+\epsilon}(\Omega)}
    + \| \lambda\Ga'(\chi\Ga)\chi_{\Gamma,t} \|_{L^{3+\epsilon}(\Gamma)}
    + \| h\Ga \|_{L^{3+\epsilon}(\Gamma)} \Big)
\end{equation}   
and we have to estimate the \rhs\ over the ``small'' time
interval $(0,S)$. Actually, the terms depending on ${\bf H}$
are controlled thanks to assumption \eqref{hypH+new}. 
We just give an estimate
for the bulk term depending on $\chi_t$, which is the most
difficult one due to worse embeddings holding in 3D. 
To do this, we use estimate \eqref{contok3} with 
$p=9+3\epsilon$ (this value is selected just for simplicity
of computation). 
%
%
Then, using assumption \eqref{hyplambda} and interpolation,
\begin{equation}\label{contox4}
  \| \lambda_b'(\chi)\chi_t \|_{L^{3+\epsilon}(\Omega)}
   \le c \big( 1 + \| \chi \|_{L^{9+3\epsilon}(\Omega)} \big)
    \| \chi_t \|_{L^{\frac{9+3\epsilon}2}(\Omega)}
   \le c \| \chi_t \|_{L^{\frac{9+3\epsilon}2}(\Omega)}.
\end{equation}   
Thus, using that
\begin{equation}\label{contox5}
  H^{\frac{5+3\epsilon}{6+2\epsilon}}(\Omega)
   \subset L^{\frac{9+3\epsilon}2}(\Omega)
\end{equation}   
and standard interpolation properties of Sobolev spaces, 
we obtain
\begin{equation}\label{contox6}
  \| \chi_t \|_{\frac{9+3\epsilon}2}
   \le c \| \chi_t \|_{V}^{\frac{5+3\epsilon}{6+2\epsilon}}
     \| \chi_t \|^{\frac{1-\epsilon}{6+2\epsilon}}\\
   \le c \Big( \| \chi_t \|_{V}^{\frac{5+3\epsilon}{6+2\epsilon}}
         t^{\frac{5+3\epsilon}{12+4\epsilon}} \Big)
    \Big( \| \chi_t \|^{\frac{1-\epsilon}{6+2\epsilon}}
         t^{\frac{1-\epsilon}{12+4\epsilon}}\Big)
        t^{-\frac12},
\end{equation}
whence, using estimate \eqref{local_est},
\begin{equation}\label{contox7}
  \| \lambda_b'(\chi)\chi_t \|_{L^{3+\epsilon}(\Omega)}
   \le c \| \chi_t \|_{\frac{9+3\epsilon}2}
   \le c t \| \chi_t \|_{V}^2 
    + c t^{-\frac{6+2\epsilon}{7+\epsilon}}.
\end{equation}   
Hence, choosing $\epsilon>0$ small enough so that the latter
exponent of $t$ is strictly larger than $-1$, we can integrate
\eqref{contox2} over $(0,S)$ to obtain 
\begin{equation}\label{st-41}
  \| \teta \|_{L^\infty(0,S;L^{3+\epsilon}(\barO,\dm))} 
   \le Q\big(\EE_0,\|\teta_0\|_{L^{3+\epsilon}(\barO,\dm))}
       ,\|\chi_0\|_{L^{\infty}(\barO,\dm)},S\big).
\end{equation}
Hence, \eqref{inpiu-2b} holds and we can apply Lemma~\ref{moser-teta}
(for $S>0$ arbitrarily small) to get \eqref{asym-teta}, 
as desired.~~\dimbox

\smallskip

\noindent%
{\bf Proof of Corollary~\ref{cor:strong_heat}.}~~%
First, note that, by \eqref{asym-teta}, 
for all $\tau>0$ we have 
$$ 
  \vert u(x,t)\vert + \vert u_\Gamma(x,t)\vert\ge c(\tau) 
     \hbox{  for a.a.~} (x,t) \in \overline\Omega\times (\tau,T),
$$
with $c(\cdot):\mathbb{R}^+\longrightarrow \mathbb{R}^+$
possibly going to~$0$ when $\tau\searrow 0$.
Hence, using \eqref{local_est}, we deduce that
\begin{equation} \label{strong_heat1}
  \int_{\tau}^T \| \teta_t(s)\|^2\ds\le c(\tau)^{-2} \tau^{-1} 
      \int_{\tau}^T s \int_\Omega\teta_t^2 u^2\,\dix\dis
   \le  Q(\mathbb{E}_0,\tau^{-1},T)
\end{equation}
as well as
\begin{equation} \label{strong_heat1ga}
  \int_{\tau}^T \| \eta_t(s)\|^2_\Gamma \ds
   \le c(\tau)^{-2}\tau^{-1} \int_{\tau}^T s \int_\Gamma\eta_t^2 u^2_\Gamma\dS\ds
   \le Q(\mathbb{E}_0,\tau^{-1},T).
\end{equation}
These relations correspond to \eqref{tetaH1}. To conclude,
we need to prove that \eqref{heat} can be interpreted in the 
strong form~\eqref{57}. Actually, this follows by reasoning
as done in the ``Estimate of the nonlinear terms''.
Namely, a comparison of terms in the bulk component of the 
heat equation gives that $\Delta u \in L^2(\tau,T;H)$.
In turn, this provides additional regularity of $\dn u$
and, as a consequence, that $\Delta\Ga u\Ga \in L^2(\tau,T;H\Ga)$. 
The details are left to the reader.~~\dimbox


\section{Long-time behavior}
\label{sec:long}

In this section we prove existence and regularity properties
of $\omega$-limit sets of weak solutions to Problem~(P). Moreover,
we also prove existence of the global attractor. 
For simplicity, we will assume that no external 
heat source is present, i.e., we will take ${\bf H}=0$. 
An asymptotically vanishing source could be treated
as well, paying the price of technical complications.
Moreover, we will restrict ourselves to the
case of singular potentials (the case of 
regular potentials being in fact simpler)
and, in particular, we will assume the reinforced
compatibility condition \eqref{specsing}. 
Finally, we will assume here the higher summability
property \eqref{hyptetaLp} on $\teta_0$
(cf.~Remark~\ref{whyLp} below for a motivation
for this choice).

To start with, we write the stationary problem associated to
Problem~(P). In order to properly state it,
we have to notice that, testing \eqref{heat} by $1$,
one gets that the quantity
\begin{equation}\label{mean}
  \mu = \mu(\teta,\chi) := \io \big( \teta + \lambda_b(\chi) \big)
   + \iga \big( \eta + \lambda\Ga(\chi\Ga) \big),
\end{equation}
representing the ``total mass'' of the (bulk+boundary)
internal energy, is conserved in the time-evolution of the system. 
Hence, once a solution trajectory evolves from some initial datum
$(\teta_0,\chi_0)$ having finite energy $\EE_0$, any point
in the $\omega$-limit set must respect the constraint
\eqref{mean}, with $\mu$ depending on the given initial datum. 

Moreover, it is apparent that the stationary version
of \eqref{heat} simply prescribes that $u$ is a constant.
Hence, we may write the steady state problem associated
to Problem~(P) as the following system:
\begin{align}\label{heatstaz}
  & \teta\ii \in (0,+\infty), ~~u\ii = -1/\teta\ii,\\
 \label{phasestaz1}
  & - \Delta \chi\ii + f(\chi\ii) - \delta \chi\ii = \lambda_b'(\chi\ii) u\ii,\\
 \label{phasestaz2}
  & - \Delta\Ga \chi\iig + f\Ga(\chi\iig) - \delta \chi\iig 
     = \lambda\Ga'(\chi\iig) u\ii - \dn \chi\ii,\\
 \label{vincolo}
  &  \io \big( \teta\ii + \lambda_b(\chi\ii) \big)
   + \iga \big( \teta\ii + \lambda\Ga(\chi\iig) \big) = \mu.
\end{align}
System \eqref{heatstaz}-\eqref{vincolo} will be named as 
Problem~(P$\iim$) in what follows. We can now state 
and prove our first result
regarding existence of nonempty $\omega$-limit sets:
\bete\label{teo:infty1}
 Let the assumptions of\/ {\rm Theorem~\ref{teor:existenceL2}} hold
 and let ${\bf H}=0$. Moreover, let\/ \eqref{specsing} 
 and\/ \eqref{hyptetaLp} hold.
 Let $(\teta,\chi)$ be the corresponding weak solution 
 to\/ {\rm Problem~(P)}. Then, as $t\nearrow\infty$,
 $(u(t),\chi(t))$ is precompact in $\calH\times \calV$.
 Moreover, any limit point $(\teta\ii,\chi\ii)$ of 
 (any subsequence of) $(\teta(t),\chi(t))$ is a solution of\/
 {\rm Problem (P$\iim$)}, where the quantity $\mu$ in\/
 \eqref{vincolo} is equal to the ``mass'' $\mu(\teta_0,\chi_0)$ 
 of the initial internal energy.
\ente
\noindent%
{\bf Proof of Theorem~\ref{teo:infty1}.}~~%
As a first step, we need to go back to the estimates performed
in Section~\ref{sec-proof-exiuni}. In particular, repeating the
energy estimate~\eqref{energy_est} and noting that now ${\bf H}=0$,
we immediately obtain 
\begin{align}\label{stunif-11}
  & \| \teta - \log \teta \|_{L^\infty(0,+\infty;L^1(\barO,\dm))} 
   + \| \nabla u \|_{L^2(0,+\infty;H)} 
   + \| \nabla\Ga u\Ga \|_{L^2(0,+\infty;H\Ga)}
   \le C,\\
 \label{stunif-12}
  & \| \chi_t \|_{L^2(0,+\infty;\calH)} 
   + \| \chi \|_{L^\infty(0,+\infty;\calV)} 
   + \| \chi \|_{L^\infty((0,+\infty)\times\barO)} 
  \le C.
\end{align}
Here and below, $C>0$ denotes suitable constants possibly
depending on $\EE_0$ but assumed in any case to be independent of 
the time variable. Instead, we will denote by $c_\mu>0$ the
constants, also independent on time, that are allowed to depend on 
the initial data only through the conserved value $\mu$.

Thanks to \eqref{poinc-u}-\eqref{poinc-uga},
we also get  
\begin{equation}\label{stunif-13}
  \| u \|_{L^2(t,t+1;\calV)} \le C.
\end{equation}
In particular, \eqref{stunif-11}-\eqref{stunif-13}
imply that 
%
%
\begin{equation}\label{inizS}
  \perogni t\ge 0,~\esiste S\in [t,t+1]:~  
   \| u(S) \|_{\calV}
   + \| \chi_t(S) \|_{\calH} 
  \le C.
\end{equation}
Then, we can go back to the ``Second estimate'' of
Subsec.~\ref{subsec:rego}. Actually, we test \eqref{heat}
by $u_t$ and the time derivative of \eqref{phase} by $\chi_t$
(hence, we do not need here the weight $t$). 
This leads to the analogue of \eqref{conto31}, namely,
\begin{align}\nonumber
  & \ddt \Big( \frac12 \| \nabla u \|^2 
   + \frac12\|\nabla_\Gamma u_\Gamma\|^2_\Gamma
   + \frac12 \| \chi_t \|^2_\calH \Big)
   + \io \teta_t^2 u^2 
   + \iga\eta_t^{2} u_\Gamma^{2}
   + \| \nabla \chi_t \|^2 
   + \|\nabla_\Gamma \chi_{\Gamma,t} \|^2_\Gamma \\
 \label{conto31u} 
  & \mbox{}~~~~~
   \le \io\lambda''_b(\chi)\chi_t^{2} u 
    + \iga \lambda''_{\Gamma}(\chiga) (\chi_{\Gamma,t})^2 u_\Gamma.
\end{align}
Then, estimating the above \rhs\ as in \eqref{conto31bis},
integrating over $(S,S+2)$, $S$ as in \eqref{inizS}, 
and using the generalized Poincar\'e inequality
\eqref{log-poincare}, we readily get
\begin{align}\label{stunif-14}
  & \| \chi_t \|_{L^\infty(1,+\infty;\calH)} 
   + \| u \|_{L^\infty(1,+\infty;\calV)} 
  \le C,\\
 \label{stunif-14b}
  & \| (\log\teta)_t \|_{L^2(t,t+1;\calH)} 
   + \| \chi_t \|_{L^2(t,t+1;\calV)} 
  \le C,~~\perogni t\in [1,+\infty).
\end{align}
Going back to \eqref{conto31u}, and having now \eqref{stunif-14} 
at our disposal, we note that the $\calV$-norm of $u$ on the \rhs\ of
\eqref{conto31bis} is controlled uniformly in time. Hence, 
we can integrate \eqref{conto31u} over the whole
$(1,+\infty)$, so that \eqref{stunif-14} is improved to
\begin{equation}\label{stunif-14x}
  \| \chi_t \|_{L^2(1,+\infty;\calV)} 
   + \| (\log\teta)_t \|_{L^2(1,+\infty;\calH)} 
   \le C.
\end{equation}
Next, we can repeat estimate \eqref{giulio11x} 
over $(S,S+2)$, obtaining
\begin{equation}\label{stunif-15}
  \| f(\chi) \|_{L^\infty(1,+\infty;H)} 
   + \| f\Ga(\chi\Ga) \|_{L^\infty(1,+\infty;H\Ga)} 
  \le C.
\end{equation}
This allows us, as before, to get also
\begin{equation}\label{stunif-16}
  \| \chi \|_{L^\infty(1,+\infty;H^2(\Omega))} 
   + \| \chi\Ga \|_{L^\infty(1,+\infty;H^2(\Gamma))} 
  \le C.
\end{equation}
Next, \eqref{stunif-11}-\eqref{stunif-12} 
and a comparison of terms in \eqref{heat} 
allow us to see that 
\begin{equation}\label{stunif-17}
  \| \teta_t \|_{L^2(0,+\infty;\calV')}
   \le C.
\end{equation}
Thanks also to the Aubin-Lions lemma, the first of
\eqref{stunif-14}, \eqref{stunif-16},
and the compact embedding $\calW\subset \calV$,
immediately give the precompactness of 
$\chi(t)$ in $\calV$. 

Applying Lemma~\ref{lemma:moser-u} over the generic
interval $(S,S+2)$, we also infer that 
\begin{equation}\label{stunif-18}
  \| u \|_{L^\infty((1,+\infty)\times \barO)} 
   \le C,
\end{equation}
which, combined with \eqref{stunif-14x},
gives
\begin{equation}\label{stunif-19}
  \| u_t \|_{L^2(1,+\infty;\calH)}
   \le C.
\end{equation}
Thanks to \eqref{stunif-18}, we can repeat the comparison
argument of Corollary~\ref{cor:sepsing}. In particular, we get
the {\sl separation property}\/ uniformly in time:
\begin{equation} \label{separunif}
  | \chi(t,x) | \le 1 - \omega \ \text{$\dm$-almost everywhere
   in $\barO$ for a.e.~$t\ge 1$,}
\end{equation}
with $\omega\in(0,1)$ independent of time.

By \eqref{stunif-19}, \eqref{stunif-14} and the 
Aubin-Lions lemma we deduce that $u(t)$ is 
precompact in $\calH$. 
With the precompactness of $(u(t),\chi(t))$ and
the above estimates at our disposal, we can
now prove that the trajectory $(\teta(t),u(t))$ 
admits a nonempty $\omega$-limit set. 

Namely, we take $\{t_n\}$ to be a diverging sequences of
time and we consider Problem~(P) over the time interval
$(0,1)$ and with initial datum $(\teta(t_n),\chi(t_n))$.
Let us call Problem~(P$_n$) this problem. Hence, once
$(\teta,\chi)$ is a weak solution to (P), 
$(\teta_n(t),\chi_n(t)):=(\teta(t_n+t),\chi(t_n+t))$, 
$t\in(0,1)$, solves~(P$_n$).

Thanks to the previous estimates and to the Aubin-Lions
lemma, there exist limit functions $\ubar,\chibar$
defined over $(0,1)\times\barO$ such that, say,
\begin{align}\label{coun11}
  & u_n \to \ubar \ \text{strongly in }
    C^0([0,1];H^{1-\epsi}(\Omega)\times H^{1-\epsi}(\Gamma)),\\
 \label{coun12}
  & \chi_n \to \chibar \ \text{strongly in }C^0([0,1];\calV),
\end{align}
for all $\epsi>0$.
As usual, all convergences are to be intended up to the 
extraction of (nonrelabelled) subsequences of $n\nearrow\infty$.
Moreover, thanks to \eqref{stunif-14x} and \eqref{stunif-19}, 
both $\ubar$ and $\chibar$ are constant in time; hence, they
coincide with the limit of (the extracted subsequence of)
$u(t_n)$ and $\chi(t_n)$, respectively.
Moreover, thanks to \eqref{stunif-11}, $\ubar$ is also constant
in space and, thanks to \eqref{stunif-18}, we have 
$|\ubar| \le C$ for some $C>0$ possibly depending on 
$\EE_0$.

Hence, passing to the limit in \eqref{phase},
we get that $(\ubar,\chibar)$ solve 
\eqref{phasestaz1}-\eqref{phasestaz2}, which can be written
in the strong form thanks to the same considerations on regularity as 
those made for the evolutionary system. In particular, the limits of the
terms $f(\chi_n)$ and $f(\chi_{n,\Gamma})$ can be identified
as before by monotonicity methods. 

To conclude the proof, we have to take the limit in the 
heat equation \eqref{heat} in order to recover \eqref{heatstaz}
and \eqref{vincolo}. To do this, we recall that, for all $t\in[0,1]$,
we have 
\begin{equation}\label{coun-13}
  \teta_n(t) = - \frac1{u_n(t)}, \quad  
   \eta_n(t) = - \frac1{(u_n(t))\Ga}, 
\end{equation}
almost everywhere in $\Omega$ and, respectively, on $\Gamma$.
Then, by \eqref{coun11}, we have  
\begin{equation}\label{coun-14}
  \teta_n(t) \to - \frac1{\ubar^{-1}}, \quad
     \eta_n(t) \to - \frac1{\ubar^{-1}},
\end{equation}
almost everywhere in $\Omega$ and, respectively, on $\Gamma$.
Let us now prove that, actually, $\ubar$ cannot be $0$.
Indeed, by \eqref{mean}, we have, for 
almost all $t\in[0,1]$,
\begin{equation}\label{coun-15}
  \ibaro \teta_n(t) \dm
   \le | \mu | 
    + | \Omega | \max_{r\in[-1,1]} |\lambda_b(r)|
    + | \Gamma | \max_{r\in[-1,1]} |\lambda\Ga(r)|
   \le c_\mu.
\end{equation}
Consequently, applying Jensen's inequality,
we deduce that there exists $c_\mu>0$ depending
only on $\mu$, $\lambda_b$ and $\lambda\Ga$ such that
$|\ubar|\ge c_\mu$. To characterize the limit 
of $\teta_n$, we need the following lemma:
\bele\label{teo:infty2}
 Under the assumptions of\/ {\rm Theorem~\ref{teo:infty1}},
 if we have in addition 
 \begin{equation}\label{teta3-11}
   (\teta_0,\eta_0) \in L^p(\Omega)\times L^p(\Gamma),
    \ \text{for some }p\ge 3,
 \end{equation}
 then it follows that
 \begin{equation}\label{teta3-12}
   \| (\teta(t),\eta(t)) \|_{L^p(\Omega)\times L^p(\Gamma)} \le 
     Q_p\big(\EE_0,\|\teta_0\|_{L^p(\Omega)},\|\eta_0\|_{L^p(\Gamma)}\big),~~
    \perogni t\ge 0,
 \end{equation}
 with $Q_p$ independent of $t\in[0,+\infty)$.
\enle
\begin{proof}
For simplicity, we just give the proof in the case when $p=3$, 
the general case following by repeating the procedure and 
applying a simple bootstrap argument.
Firstly, let us observe that, by \eqref{st-41} with $\epsilon=0$, we have
\begin{equation}\label{teta3-13}
  \| (\teta(t),\eta(t)) \|_{L^3(\Omega)\times L^3(\Gamma)} 
   \le Q\big(\EE_0,\|\teta_0\|_{L^3(\Omega)},\|\eta_0\|_{L^3(\Gamma)}\big)
   ~~\perogni t\in [0,1].
\end{equation}
Hence, we can test \eqref{heat} by $\teta^2$. Using that now $\chi$ is
uniformly bounded, we obtain (compare with \eqref{contox2}),
for some $\kappa>0$ (which may change from line to line in the 
computations below),
\begin{equation}\label{teta3-14}
  \ddt \| \teta \|_{L^{3}(\barO,\dm)}^3
   + \kappa \| \nabla \teta^{1/2} \|^2
   + \kappa \| \nabla\Ga \eta^{1/2} \|^2\Ga
  \le c \ibaro | \chi_t \teta^2 | \dm.
\end{equation}   
Actually we can notice that, if \eqref{teta3-11} holds, then
the solution is regular enough so that $\eta$ is in fact the trace 
of $\teta$ in this case.

Noting that, due to \eqref{mean},
\begin{equation}\label{teta3-15}
  \| \teta^{1/2}(t) \|_{L^2(\barO,\dm)}^2 
   = \ibaro \teta(t) \dm \le c_\mu
   ~~\perogni t\in[0,+\infty),
\end{equation}
adding \eqref{teta3-15} to \eqref{teta3-14}, and using
standard Sobolev's embeddings, we infer
\begin{align}\nonumber
  \ddt \| \teta \|_{L^3(\barO,\dm)}^3
   + \kappa \| \teta \|_{L^3(\barO,\dm)}
   & \le c \| \chi_t \|_{L^3(\barO,\dm)}
    \| \teta^{3/2} \|_{L^2(\barO,\dm)}
    \| \teta^{1/2} \|_{L^6(\barO,\dm)}
    + c_\mu \\
 \label{teta3-16}
   & \le c \| \chi_t \|_{L^3(\barO,\dm)}^2
    \| \teta \|_{L^3(\barO,\dm)}^3
    + \frac\kappa2 \| \teta \|_{L^3(\barO,\dm)} 
    + c_\mu.
\end{align}
Thus, setting
\begin{equation}\label{teta3-17}
  y(t) := \max\Big\{ 1, \| \teta(t) \|_{L^3(\barO,\dm)}^2 \Big\},
   \ \ m(t):= \| \chi_t(t) \|_{L^3(\barO,\dm)}^2,
\end{equation}
it is clear that \eqref{teta3-16} can be interpreted as the
differential inequality
\begin{equation}\label{teta3-18}
  y'(t) + \kappa \le 
    c m(t) y(t) + c_\mu y(t)^{-1/2},
\end{equation}
so that the comparison principle, together with 
\eqref{teta3-13}, \eqref{stunif-14x} and Gronwall's lemma,
readily imply that $y$ is uniformly bounded
for $t\in[1,+\infty)$, as desired.
Indeed, in the set of times such that
$y(t)\le 4c_\mu^2/\kappa^2$ there is nothing to prove,
while if $y> 4c_\mu^2/\kappa^2$ \eqref{teta3-18}
takes the form $y' \le c m y$ and we can use 
Gronwall's lemma since $m$ is globally summable 
thanks to \eqref{stunif-14x}.
\end{proof}
\noindent%
\beos\label{whyLp}
 It is worth noticing that the property stated in the 
 Lemma is likely to be false if $p\in(1,3)$. Indeed,
 in that case the gradient terms in \eqref{teta3-14} 
 seem to give no help and, consequently, the differential
 inequality corresponding to \eqref{teta3-18} loses its
 dissipative character (in other words, it provides
 a control on the $L^p$-norm of $\teta$ which is not
 uniform in time). This is the reason which led us
 to assume \eqref{hyptetaLp} in this section. 
 On the other hand, thanks to the conservation of 
 $\mu$ (cf.~\eqref{mean}), we have in any case a uniform
 in time $L^1$-control on $\teta$. However,
 this seems not enough in order to characterize
 properly the limit of $\teta_n$. Indeed, 
 it may happen that
 $$ 
   \teta_n = \teta_{n,1} + \teta_{n,2}, \ \ 
    \text{so that } u_n = -\frac1{\teta_{n,1} + \teta_{n,2}},
 $$
 where $\teta_{n,1}$ is a ``good'' function which converges, say,
 uniformly, to the constant $-1/\ubar$. Instead, the sequence
 $\{\teta_{n,2}\}$ is also bounded in $L^1$ but it may
 converge to some singular measure (e.g., a Dirac mass, 
 meaning that the support of 
 $\teta_{n,2}$ becomes small and $\teta_{n,2}$ becomes large in
 that set). However, when one computes $u_n$, the contribution of
 $\teta_{n,2}$ is negligible for large $n$ because $\teta_{n,2}$ 
 influences the value of $u_n$ only in a set that has asymptotically 
 measure $0$. Hence, the limit of $\teta_n$ is in this case 
 a measure whose regular (absolutely continuous)
 part coincides with $-1/\ubar$, but which may also
 have a singular component.
\eddos
\noindent%
%
%
%
%
%
Thanks to Lemma~\ref{teo:infty2}, using pointwise convergence
(cf.~\eqref{coun-14}) and Lebesgue's theorem, we obtain
\begin{equation}\label{coun-21}
  \teta_n \to - \frac1{\ubar^{-1}} \ \text{strongly in }
    L^p(0,1;L^p(\barO;\dm))~~\perogni p\in[1,3).
\end{equation}
In particular, we get the second \eqref{heatstaz} 
%
%
\eqref{vincolo}, which concludes the proof of 
Theorem~\ref{teo:infty1}.~~\dimbox

\medskip 

\noindent%
Finally, we come to the problem of the existence of the 
global attractor for Problem~(P). To address this issue,
we first introduce the natural {\sl phase space}\/
for the dynamical system associated to Problem~(P), by
setting
\begin{equation}\label{attra11}
  \calX_\mu:=\left\{ (\teta,\chi) \in 
     (L^{3+\epsilon}(\Omega)\times L^{3+\epsilon}(\Gamma)) \times \calV:
     \begin{cases}
	\teta>0~\text{$\dm$-a.e.~in $\barO$},\ \log\teta\in L^1(\barO;\dm)\\
	F(\chi) \in L^1(\Omega),\ F\Ga(\chi\Ga) \in L^1(\Gamma),\\
        \io \big( \teta + \lambda_b(\chi) \big)
	    + \iga \big( \eta + \lambda\Ga(\chi\Ga) \big) = \mu
      \end{cases}   \right\},
\end{equation}
where $\epsilon>0$ is given and
the conserved value $\mu$ is also prescribed. It is not difficult
to prove (see, e.g., \cite{RS} for details) that the above space,
endowed with the distance
\begin{align}\nonumber
  & \diX((\teta_1,\chi_1),(\teta_2,\chi_2))
   := \| \teta_1 - \teta_2 \|_{L^{3+\epsilon}(\barO,\dm)}
    + \| \log \teta_1 - \log \teta_2 \|_{L^1(\barO,\dm)}\\
 \label{defdX}
  & \mbox{}~~~~~
    + \| \chi_1 - \chi_2 \|_{\calV}
    + \| F(\chi_1) - F(\chi_2) \|_{L^1(\Omega)}
    + \| F\Ga(\chi_{1,\Gamma}) - F\Ga(\chi_{2,\Gamma}) \|_{L^1(\Gamma)}
\end{align}
acquires a complete metric structure.

We will name Problem~(P$_\mu$) the version 
of Problem~(P) where the ``mass'' $\mu$ of the initial
internal energy is assigned.
\beos\label{muammi}
 In order for the phase space $\calX_\mu$ not to be empty, we 
 implicitly assume that  
 \begin{equation}\label{minimu}
   \mu > | \Omega | \min_{r\in[-1,1]} \lambda_b(r)
    + | \Gamma | \min_{r\in[-1,1]} \lambda\Ga(r).
 \end{equation}
 Actually, if \eqref{minimu} is not satisfied, then
 from \eqref{mean} we get that $\teta$ needs to be nonpositive,
 and the problem becomes inconsistent. However, we will
 see just below that a stronger condition may be needed. 
\eddos
%
%
\bete\label{teo:attra}
 Let the assumptions of\/ {\rm Theorem~\ref{teo:infty1}} 
 hold. Let also assume that, either 
 \begin{equation}\label{chiavemu}
   \mu > | \Omega | \max_{r\in[-1,1]} \lambda_b(r)
     + | \Gamma | \max_{r\in[-1,1]} \lambda\Ga(r), 
 \end{equation}
 or
 \begin{equation}\label{segnola}
   \liminf_{|r|\nearrow1} \lambda'(r)\sign r > 0.
 \end{equation}
 Then, the dynamical process associated 
 to~{\rm Problem~(P$_\mu$)} admits the global
 attractor $\calA_\mu$. Moreover, there exist constants
 $K>0$, $\omega\in(0,1)$ and $\alpha\in(0,1)$ depending only
 on $\mu$ such that, for any $(\teta,\chi)\in\calA_\mu$,
 \begin{align}\label{regA1}
   & \| \teta \|_{\calV} 
    + \| \chi \|_{\calW} \le K,\\
  \label{regA2}
   & \alpha \le \teta \le \alpha^{-1}
    ~~\text{$\dm$-almost everywhere in $\barO$},\\
  \label{regA3}
   & - 1 + \omega \le \chi \le 1 - \omega
    ~~\text{$\dm$-almost everywhere in $\barO$}.   
 \end{align}
\ente
\begin{proof}
Due to the singular and degenerate character of equation
\eqref{heat}, a direct proof of a dissipative estimate 
(e.g., of the existence of an absorbing set) appears to
be out of reach. On the other hand, as in \cite{Sc09}, we 
can take advantage of the fact that our system admits a
coercive Liapounov functional (the energy $\calE$). 
Thanks to this property, the existence of the global attractor
follows by proving the following conditions:
\begin{itemize}
 \item[(a)] Solution trajectories are precompact with respect to
  the metrics $\diX$ of the phase space $\calX_\mu$;
 \item[(b)] The set of stationary states is bounded in the energy space
  {\sl independently} of the magnitude of the initial data.
\end{itemize}
Hence, let us start with the proof of (a). By
Theorem~\ref{teo:infty1}, it is clear that for any weak solution
there exists $C$ depending on $\EE_0$ and on 
the $L^{3+\epsilon}$-norm of $\teta_0$ such that
\begin{equation}\label{un11}
  \| u(t) \|_{\calV}
   + \| \chi(t) \|_{\calW}
   \le C~~\perogni t\in[1,+\infty).
\end{equation}
Moreover, thanks also to Lemma~\ref{teo:infty2}, 
we have 
\begin{equation}\label{un12}
  0 < \alpha \le \teta(t,x) \le \alpha^{-1}
   \ \ \text{$\dm$-a.e.~in}~\barO~~
   \text{and for all } t\in[1,+\infty),
\end{equation}
with $\alpha\in(0,1)$ depending on the same quantities
as the above constant $C$. Coupling \eqref{un11} with \eqref{un12},
we readily get
\begin{equation}\label{un13}
  \| \teta(t) \|_{\calV}
   \le C~~\perogni t\in[1,+\infty).
\end{equation}
In addition to that, the {\sl uniform}\/
separation property \eqref{separunif} holds. Hence, both
$\teta$ and $\chi$ are eventually separated from 
the singularities.
Using \eqref{un11}-\eqref{un13} and \eqref{separunif} it 
is a standard matter to prove that the
trajectory is precompact with respect to 
the metric structure of~$\calX_\mu$.

\smallskip

To conclude the proof, let us demonstrate~(b). This is more
delicate, since we need to prove that the set
of stationary states (i.e., of the solutions $(\teta\ii,\chi\ii)$ 
to Problem~(P$_{\infty,\mu}$)) is bounded in a way that
depends only on $\mu$. In other words, all elements of 
the $\omega$-limit of a given solution trajectory are bounded
in a way that depends on the initial datum only through
the conserved quantity $\mu$. 
Actually, we have already seen 
(see \eqref{coun-15}) that $|\teta\ii|\le c_\mu$.
To prove that the same property holds for $u\ii$,
we first consider the situation when \eqref{chiavemu}
holds. In this case, combining \eqref{chiavemu} with 
\eqref{vincolo}, $|u\ii| \le c_\mu$ follows 
immediately. 

Hence, the terms $\lambda_b'(\chi\ii)u\ii$ and
$\lambda\Ga'(\chi_{\infty,\Gamma})u\iig$ on the \rhs s
of \eqref{phasestaz1} and, respectively, \eqref{phasestaz2}
are uniformly bounded in a way depending only on $\mu$.
Then, by the maximum principle it immediately follows
that there exists $\omega_\mu \in (0,1)$ such that
$-1+\omega_\mu \le \chi\ii(t,x) \le 1-\omega_\mu$
$\dm$-almost everywhere in $\barO$. 
It is immediate to check that the
same property holds if \eqref{segnola} is satisfied
in place of \eqref{chiavemu}.

Standard elliptic regularity estimates applied to
\eqref{phasestaz1}-\eqref{phasestaz2} then give, say,
\begin{equation}\label{un14}
  \| \chi\ii \|_{\calW} \le c_\mu
\end{equation}
for every stationary solution $(\teta\ii,\chi\ii)$
(actually, much more is true since $\chi\ii$ is separated
from the singularities and is, consequently, a classical solution
to the elliptic system \eqref{phasestaz1}-\eqref{phasestaz2}).
Hence, (b) holds true and the existence of the attractor
with the properties \eqref{regA1}-\eqref{regA3}
follows from classical results (see, e.g., \cite{Hale}).
\end{proof}


\begin{thebibliography}{99}











\bibitem{At} 
 H.~Attouch,
 ``Variational Convergence for Functions and Operators'',
 Applicable Mathematics Series,
 Pitman (Advanced Publishing Program), Boston, MA,
 1984.

\bibitem{Ba} 
 V.~Barbu,
 ``Nonlinear Semigroups and Differential Equations in Banach Spaces'',
 Noord\-hoff,
 Leyden,
 1976.

\bibitem{bonf_vaz2010}
M.~Bonforte, and J.~L.~ V{\'a}zquez,
 \newblock {\sl Positivity, local smoothing, and {H}arnack inequalities for
              very fast diffusion equations},
  \newblock  Adv. Math., {\bf 223}, 2010, 2, 529--578.




\bibitem{Br} 
 H.~Br\'ezis,
 ``Op\'erateurs Maximaux Monotones et Semi-groupes de Contractions
   dans les Espaces de Hilbert'',
 North-Holland Math.\ Studies {\bf 5},
 North-Holland,
 Amsterdam,
 1973.



 \bibitem{BG}
 F.~Brezzi and G.~Gilardi,
 FEM Mathematics,
 in Finite Element Handbook (H.~Kardestuncer Ed.), Part I: Chapt.~1:
 Functional Analysis, 1.1--1.5; Chapt.~2: Functional Spaces, 2.1--2.11;
 Chapt.~3: Partial Differential Equations, 3.1--3.6,
 McGraw-Hill Book Co., New York, 1987.

\bibitem{Ca}
 G.~Caginalp,
 {\sl An analysis of a phase field model of a free boundary},
 Arch.\ Rational Mech.\ Anal.,
 {\bf 92} (1986),
 205--245.

\bibitem{CGGM}
 C.~Cavaterra, C.G.~Gal, M.~Grasselli, and A.~Miranville, 
  {\sl Phase-field systems with nonlinear coupling and dynamic boundary conditions},
 Nonlinear Anal., {\bf 72} (2010), 
 2375--2399.

\bibitem{CGM} 
 L.~Cherfils, S.~Gatti, and A.~Miranville, 
 {\sl Existence of global solutions to the Caginalp phase-field system with
    dynamic boundary conditions and singular potentials},
 J.~Math.\ Anal.\ Appl., 
 {\bf 343} (2008), 
 557--566 
 [Corrigendum, J.~Math.\ Anal.\ Appl., {\bf 348} (2008), 1029--1030].

\bibitem{CM} 
 L.~Cherfils and A.~Miranville, 
 {\sl On the Caginalp system with dynamic boundary conditions 
  and singular potentials},
 Appl.\ Math.,
 {\bf 54} (2009), 
 89--115.

 \bibitem{CLau}
 P.~Colli and Ph.~Lauren\c cot,
 {\sl Weak solutions to the Penrose-Fife phase field model for a class of admissible heat flux laws},
 Phys.~D, {\bf 111} (1998), 
 311--334.

\bibitem{CoGM} 
 M.~Conti, S.~Gatti, and A.~Miranville,
 {\sl Asymptotic behavior of the Caginalp phase-field system with coupled 
   dynamic boundary conditions},
 Discrete Contin.\ Dyn.\ Syst.~S, 
 {\bf 5} (2012), 485--505.

\bibitem{FS}
 E.~Feireisl and G.~Schimperna,
 {\sl Large time behaviour of solutions to Penrose-Fife phase change models},
 Math.\ Methods Appl.\ Sci., 
 {\bf 28} (2005), 
 2117--2132.

\bibitem{FiMD1}
 H.P.~Fischer, P.~Maass, and W.~Dieterich, 
 {\sl Novel surface modes in spinodal decomposition},
 Phys.\ Rev.\ Letters, {\bf 79} (1997), 
 893--896.

\bibitem{FiMD2}
 H.P.~Fischer, P.~Maass, and W.~Dieterich,
 {\sl Diverging time and length scales of spinodal decomposition modes in thin flows},
 Europhys.\ Letters, {\bf 42} (1998), 
 49--54.

\bibitem{FRDGMMR} 
 H.P.~Fischer, J.~Reinhard, W.~Dieterich, J.-F.~Gouyet, P.~Maass, A.~Majhofer, 
   and D.~Reinel, 
 {\sl Time-dependent density functional theory and the kinetics of lattice gas 
  systems in contact with a wall},
 J.\ Chem.\ Phys., 
 {\bf 108} (1998), 
 3028--3037.

\bibitem{GGbis}
 C.G.~Gal and M.~Grasselli, 
 {\sl On the asymptotic behavior of the Caginalp system with dynamic 
  boundary conditions}, 
 Commun.\ Pure Appl.\ Anal.,
 {\bf 8} (2009), 
 689--710.


\bibitem{GMS1} 
 G.~Gilardi, A.~Miranville, and G.~Schimperna,
 {\sl On the Cahn-Hilliard equation with irregular potentials and dynamic
 boundary conditions},
 Commun.\ Pure Appl.\ Anal., 
 {\bf 8} (2009), 881--912.

 \bibitem{GMS2}
  G.~Gilardi, A.~Miranville, and G.~Schimperna,
  {\sl Long time behavior of the Cahn-Hilliard equation with irregular potentials 
    and dynamic boundary conditions},
  Chin.\ Ann.\ Math.\ Ser.~B, 
 {\bf 31} (2010), 
 679--712.

\bibitem{GoMS} 
 G.R.~Goldstein, A.~Miranville, and G.~Schimperna,
 {\sl A Cahn-Hilliard model in a domain with non-permeable walls},
 Phys.~D, 
 {\bf 240} (2011), 
 754--766.

\bibitem{gms}
 M.~Grasselli, A.~Miranville, and G. Schimperna, 
 {\sl The Caginalp phase-field system with coupled
    dynamic boundary conditions and singular potentials}, 
 Discrete Contin. Dyn.\ Syst.,
 {\bf 28} (2010), 
 67--98.

\bibitem{GPS1}
 M.~Grasselli, H.~Petzeltov\'a, and G.~Schimperna, 
 {\sl Long time behavior of solutions to the Caginalp system with singular potential},
 Z.~Anal.\ Anwend., 
 {\bf 25} (2006), 
 51--72. 
 
 
\bibitem{Hale}
 J.K.~Hale, ``Asymptotic Behavior of Dissipative Systems'',
 Mathematical Surveys and Monographs, 25. 
 American Mathematical Society, 
 Providence, RI, 1988.
 

\bibitem{HSZ}
 W.~Horn, J.~Sprekels, and S.~Zheng,
 {\sl Global existence of smooth solutions to the Penrose-Fife model 
  for Ising ferromagnets},
 Adv.\ Math.\ Sci.\ Appl., 
 {\bf 6} (1996), 
 227--241.

\bibitem{Ito-Ken-Kubo03}
 A.~Ito, N.~Kenmochi, and M.~Kubo,
 {\sl Non-isothermal phase transition models with Neumann boundary
   conditions},
 Nonlinear Anal., 
 {\bf 53}, (2003), 
 977--996.

\bibitem{IK}
 A.~Ito and N.~Kenmochi,
 {\sl Inertial set for a phase transition model of Penrose-Fife type},
 Adv.\ Math.\ Sci.\ Appl., 
 {\bf 10} (2000), 353--374.

\bibitem{IKN}
 A.~Ito, N.~Kenmochi, and M.~Niezg\'odka,
 {\sl Phase separation model of Penrose-Fife type with Signorini boundary condition},
 Adv.\ Math.\ Sci.\ Appl., 
 {\bf 17} (2007), 
 337--356.

\bibitem{KEMRSBD}
 R.~Kenzler, F.~Eurich, P.~Maass, B.~Rinn, J.~Schropp, E.~Bohl, and W.~Dieterich, 
 {\sl Phase separation in confined geometries: solving 
   the Cahn-Hilliard equation with generic boundary conditions}, 
 Comput.\ Phys.\ Commun., 
 {\bf 133} (2001),
 139--157.

\bibitem{Lau}
 Ph.~Lauren\c cot,
 {\sl Weak solutions to a Penrose-Fife model for phase transitions},
 Adv.\ Math.\ Sci.\ Appl., 
 {\bf 5} (1995), 117--138.

\bibitem{Li} 
 J.-L.~Lions,
 ``Quelques M\'ethodes de R\'esolution des Probl\`emes aux Limites
    non Lin\'eaires'' (French),
 Dunod, Gauthier-Villars,
 Paris, 1969.


\bibitem{MiranvilleRavello}
 A. Miranville, 
 ``Some mathematical models in phase transition'',
 Lecture Notes, Ravello, 2009.

\bibitem{MZ}
A.~Miranville and S.~Zelik, 
 {\sl Robust exponential attractors for Cahn-Hilliard type equations with 
  singular potentials},
 Math.\ Methods Appl.\ Sci., {\bf 27} (2004), 
 545--582. 

\bibitem{MZCH}
 A.~Miranville and S.~Zelik, 
 {\sl The Cahn-Hilliard equation with singular potentials
  and dynamic boundary conditions},
 Discrete Contin.\ Dyn.\ Syst., 
 {\bf 28} (2010), 275--310.

\bibitem{MuRo}
 D.~Mugnolo and S.~Romanelli,
 {\sl Dirichlet forms for general Wentzell boundary conditions, analytic semigroups,
  and cosine operator functions},
 Electron.\ J.~Differential Equations
 {\bf 2006}, No.~118, 20 pp.~(electronic).

\bibitem{PenroseFife}
 O.~Penrose and P.C.~Fife,
 {\sl Thermodynamically consistent models of phase-field type for the
   kinetics of phase transitions},
 Phys.~D, {\bf 43} (1990), 
 44--62.

\bibitem{RS} 
 E.~Rocca and G.~Schimperna,
 {\sl Universal attractor for some singular phase transition systems},
 Phys.~D,
 {\bf 192} (2004), 
 279--307.

\bibitem{RSMed}
 E.~Rocca and G.~Schimperna,
 {\sl Universal attractor for a Penrose-Fife system with special heat flux law},
 Mediterr.\ J.~Math., 
 {\bf 1} (2004), 109--121.

\bibitem{SaV}
 G.~Savar\'e and A.~Visintin, 
 {\sl Variational convergence of nonlinear diffusion equations: applications to 
   concentrated capacity problems with change of phase},
 Atti Accad.\ Naz.\ Lincei Cl.\ Sci.\ Fis.\ Mat.\ Natur.\ Rend.\ Lincei (9) Mat.\ Appl.,
 {\bf 8} (1997), 
 49--89.
 
\bibitem{Sc99}
 G.~Schimperna,
 {\sl Weak solution to a phase-field transmission problem in a concentrated capacity},
 Math.\ Methods Appl.\ Sci.,
 {\bf 22} (1999), 
 1235--1254.

\bibitem{Sc09}
 G.~Schimperna,
 {\sl Global and exponential attractors for the Penrose-Fife system},
 Math.\ Models Methods Appl.\ Sci.,
 {\bf 19} (2009), 969--991.

\bibitem{ssz12}
 G.~Schimperna, A.~Segatti, and S.~Zelik,
 {\sl Asymptotic uniform boundedness of energy solutions to the Penrose-Fife model},
 J.~Evol.\ Equ., {\bf 12} (2012), 863--890.

 \bibitem{ssz13}
G.~Schimperna, A.~Segatti, and S.~Zelik,
{\sl On a singular heat equation with dynamic boundary conditions},
 paper in preparation (2013).

%

\bibitem{vaz2006}
 J.L.~V{\'a}zquez,
 Smoothing and decay estimates for nonlinear diffusion equations,
 Oxford Lecture Series in Mathematics and its Applications,
 {\bf 33}, Oxford University Press, Oxford, 2006.
 

\end{thebibliography}
\end{document}